\documentclass[twocolumn]{autart}    

\usepackage{graphicx}          
\usepackage{indentfirst}
\usepackage{amsmath}
\usepackage{amssymb} 
\usepackage[compress]{cite}

\usepackage{bm}
\usepackage{subfigure}
\usepackage{enumerate}

\usepackage{graphicx}
\usepackage{epstopdf}
\usepackage{bbm}
\usepackage{bbold}

\usepackage{amsmath}

\usepackage{lineno}
\usepackage{multirow}
\usepackage{amsfonts}
\usepackage{textcomp}
\usepackage{stfloats}

\usepackage{xcolor}
\usepackage{color}
\usepackage{threeparttable}
\usepackage[bookmarksnumbered, colorlinks, citecolor=blue!50!blue!100, linkcolor=blue]{hyperref}

\usepackage{hyperref}
\hypersetup{colorlinks=true, linkcolor=blue, breaklinks=true, urlcolor=blue, citecolor=blue}

\usepackage{float}
\usepackage{indentfirst}
\usepackage{mathrsfs}

\usepackage{algorithm}
\usepackage{algorithmicx}

\usepackage{tikz}
\usepackage{centernot}
\usepackage{version,xspace}
\usepackage{environ}
\usepackage{blkarray}
\usetikzlibrary{automata,positioning,arrows,through}

\usepackage{appendix}

	\newtheorem{defi}{\textbf{Definition}}[section]
	\newtheorem{thom}{\textbf{Theorem}}[section]
	\newtheorem{asp}{\textbf{Assumption}}[section]
	\newtheorem{rek}{\textbf{Remark}}[section]
	
	\newtheorem{lema}{\textbf{Lemma}}[section]
	
	\newtheorem{pbm}{\textbf{Problem}}[section]

\newcommand{\defiref}[1]{Definition \ref{#1}}
\newcommand{\thomref}[1]{Theorem~\ref{#1}}
\newcommand{\aspref}[1]{Assumption~\ref{#1}}

\newcommand{\corref}[1]{Corollary \ref{#1}}
\newcommand{\lemaref}[1]{Lemma \ref{#1}}

\allowdisplaybreaks[4]
\begin{document}
\begin{frontmatter}

\title{Matching-Based Capture Strategies for 3D
	Heterogeneous Multiplayer Reach-Avoid Differential Games\thanksref{footnoteinfo}} 

\thanks[footnoteinfo]{The work of R. Yan, Z. Shi, and Y. Zhong was supported in part by the National Natural Science Foundation of China under Grant 61374034,
	and in part by China Scholarship Council. This work of X. Duan and
	F. Bullo was supported in part by the US Air Force Office of Scientific
	Research under award FA9550-15-1-0138.}

\author[Beijing]{Rui Yan}\ead{yr15@mails.tsinghua.edu.cn},    
\author[UCSB]{Xiaoming Duan}\ead{xiaomingduan@ucsb.edu},  
\author[Beijing]{Zongying Shi}\ead{szy@mail.tsinghua.edu.cn}, 
\author[Beijing]{Yisheng Zhong}\ead{zys-dau@mail.tsinghua.edu.cn}, 
\author[UCSB]{Francesco Bullo}\ead{bullo@ucsb.edu} 

\address[Beijing]{Department of Automation, Tsinghua University, Beijing 100084, China}
\address[UCSB]{Center of Control, Dynamical Systems and Computation, University of California, Santa Barbara, USA}
          

\begin{abstract}                          
\hspace{0.2 in}  This paper studies a 3D multiplayer reach-avoid differential game with a goal region and a play region. Multiple pursuers defend the goal region by
consecutively capturing multiple evaders in the play region. The players
have heterogeneous moving speeds and the pursuers have heterogeneous
capture radii. Since this game is hard to analyze directly, we decompose the whole game as many subgames which involve multiple pursuers and only one evader. Then, these subgames are used as a building block for the pursuer-evader matching. First, for multiple pursuers and one evader, we introduce an evasion space (ES) method
characterized by a potential function to construct a guaranteed pursuer
winning strategy. Then, based on this strategy, we develop conditions to
determine whether a pursuit team can guard the goal region against one
evader. It is shown that in 3D, if a pursuit team is able to defend the
goal region against an evader, then at most three pursuers in the team
are necessarily needed. We also compute the value function of the
Hamilton-Jacobi-Isaacs (HJI) equation for a special subgame of degree. To capture the maximum number of
evaders in the open-loop sense, we formulate a maximum bipartite matching problem with conflict
graph (MBMC). We show that the MBMC is NP-hard and design a
polynomial-time constant-factor approximation algorithm to solve
it. Finally, we propose a receding horizon strategy for the pursuit team
where in each horizon an MBMC is solved and the strategies of the pursuers are given. We also extend our results to the case of a
bounded convex play region where the evaders escape through an exit. Two
numerical examples are provided to demonstrate the obtained results.

\end{abstract}

\end{frontmatter}

\section{Introduction}
\emph{Problem description and motivation:} Consider a multiplayer reach-avoid differential game in 3D space where a plane divides the space into two disjoint regions, i.e., the play region and the goal region.
An evasion team with multiple evaders of different speeds, initially lying in the play region, tries to send as many its team members as possible into the goal region. Meanwhile, a pursuit team with multiple pursuers of different speeds and capture radii, initially spreading over the space, aims to guard the goal region by capturing the evaders. From a different point of view, this is also equivalent to a game where the evaders try to escape from the play region and avoid adversaries and dynamic obstacles formulated as a pursuit team. As an extension, we also consider a game played in a bounded convex region in 3D with a planar exit. This game is hard to solve directly because of the high dimension and multiple stages \cite{MC-ZZ-CJT:17,PA-SDB-FB:14d,SL-ZZ-CJT-KH:13,RY-ZS-YZ:20}. This paper provides a matching-based capture strategy for the pursuit team to capture a good number of evaders, partially dealing with this game.   

\quad This problem is motivated by robotic applications, including the robot competition, dynamic collision avoidance and region surveillance \cite{FB-EF-MP-KS-SLS:10k,DWO-PTK-ARG:16,XD-MG-FB:17o}. For example, in region protection games, multiple pursuers are used to intercept multiple adversarial
intruders. In collision avoidance and path planning, how a group of vehicles can get into some target set or escape from a bounded region through an exit, while avoiding dangerous situations, such as collisions with moving obstacles.

\emph{Literature review:} The problem in this paper is related to games such as lifeline games, two-target differential games, reach-avoid games and target guarding differential games.

\quad The two-player lifeline games were introduced by Isaacs in \cite{RI:65} and then a two-pursuer-one-evader planar case in a square domain was revisited by \cite{RY-ZS-YZ:19}. By fixing the evaders' speeds and formulating them as demands, \cite{SDB-SLS-FB:08v} designed a service policy for a pursuer and derived the conditions for its stability based on system parameters. Recently, the task assignment for multiple pursuers and evaders in convex planar domains has also been studied in \cite{RY-ZS-YZ:20} by computing analytical barriers. As for two-target differential games, the first quantitative and qualitative results appeared in \cite{AB-FG-GL:69}, where each of two players has her own target toward which she wishes to steer the system state before the other. In \cite{GJO-JVB:74,MP-WMG:80,WMG-MP:81}, several variations of two-player games were considered such as role determination, complex dynamics and targets of different shapes.

\quad Reach-avoid differential games were first discussed in \cite{IMM-AMB-CJT:05,KM-JL:11,ZZ-RT-HH-CJT:12}, and then extensive studies including many variations and practical applications appeared \cite{HH-JD-WZ-CJT:15}. The current method for these games involves solving a Hamilton-Jacobi-Isaacs (HJI) equation, which suffers from the curse of dimensionality, so various techniques have been proposed, including approximation function \cite{VSC-DP:19}, system decomposition \cite{MC-SLH-MSV-SB-CJT:18}, and cone programming \cite{JL-MC-BL-MP:18}. Most of these works focus on approximation, two-player or open loop games  \cite{SL-ZZ-CJT-KH:13} due to the exponential growth of the computation as the size of the states increases. For multiplayer cases, \cite{MC-ZZ-CJT:17} greatly reduced the computation burden by creating a number of straight lines in 2D to output matching pairs.

\quad The two-player target guarding differential games were also first studied by Isaacs in \cite{RI:65}, and revisited by \cite{JM-SRM-RHV-BB:19} with the goal of real-time implementation. Recently, multiplayer cases have received increasing attention from algorithms to game setups. For example, the authors in \cite{DS-DAP:19} used the swarming behavior of male mosquitoes to design the motion strategy for multiple guardians against a fast intruder in area protection. Multiplayer scenarios of special setup were explored in \cite{RY-ZS-YZ:17-2,DS-VK:18,PA-SDB-FB:14d}, where \cite{RY-ZS-YZ:17-2} and \cite{DS-VK:18} restrict the motion of pursuers on the boundary of target set, and \cite{PA-SDB-FB:14d} studied the escape from a circular disk.

\quad To the best of our knowledge, there are only a few works concerning multiplayer reach-avoid differential games. Chen \emph{et al.} in \cite{MC-ZZ-CJT:17} considered a multiplayer
reach-avoid differential game where no cooperation except the matching among pursuers exists and thus provided a suboptimal solution. In \cite{RY-ZS-YZ:20}, the authors focused on the case of zero capture radius and homogeneous players, and they solved the task assignment problem by a 0-1 integer programming without complexity analysis and polynomial-time algorithms. In \cite{DS-VK:18,PA-SDB-FB:14d,RY-ZS-YZ:17}, the authors limited their attention to the case of zero capture radius and non-fully competitive strategies. Moreover, it is worth noting that it is desirable to develop an analytical and efficient method to analyze multiplayer pursuit-evasion games with heterogeneous capture radii, as discussed in \cite{SYH-TS:17,EG-DWC-MP:19,MC-ZZ-CJT:17,WS-PT-AJY:19}.

\quad Another focus of this work is on the \emph{constrained matching problems} \cite{AI-MR-SLT:78} or \emph{maximum matching problems with conflicts} (MMPC). In the MMPC, conflicts between edges are represented by an undirected graph, i.e., the \emph{conflict graph}. Every vertex of the conflict graph corresponds to an edge of the original graph, and every edge corresponds to a binary conflict. The MMPC is equivalent to finding the maximum matching such that at most one edge in each conflicting edge pair is selected. Conflict graphs have been considered in many combinatorial optimization problems such as knapsack problem \cite{UP-JS:17} and minimum spanning tree problems \cite{RZ-SNK-APP:11}.

\quad The first work introducing conflicts into matching problems was presented by  \cite{AI-MR-SLT:78}, where a constrained bipartite matching problem was considered. Then, Thomas \cite[Chapter 4]{DJT:16} revisited the constrained matching problem, and summarized the complexity of different variations based on fixed and variable parameters. An important work addressing the MMPC was by \cite{AD-UP-JS-GJW:11}, where they established primary complexity results for several variations. Then, these results were extended later on in \cite{TO-RZ-APP:13} where the authors proposed additional complexity results, identified special polynomially solvable cases, and also designed several heuristic algorithms. However, all current results cannot perfectly fit our problems.

\emph{Contributions:} In this paper, we study the cooperative strategies for multiple pursuers to guard a 3D region against multiple evaders.  Concretely, we propose an Evasion Space (ES) method characterized by a potential function and design a receding horizon capture strategy for the pursuit team to capture a good number of evaders. Compared with \cite{DS-VK:18,MC-ZZ-CJT:17,RY-ZS-YZ:20,DS-DAP:19,PA-SDB-FB:14d}, we consider more practical cases when the pursuers have different capture radii, can closely cooperate with each other, and have no knowledge of the strategies of the evaders. Since the evaders can take any strategies which are unavailable to the pursuers, a robust feedback pursuit strategy is needed. The existing techniques cannot be applied directly to these cases and we introduce new methods to solve them. The main contributions are as follows.
\begin{enumerate}[i)]
	\item A method is presented to study heterogeneous multiplayer reach-avoid differential games in 3D where the pursuers try to defend a goal region against evaders. Heterogeneity refers to different speeds for players and different capture radii for pursuers. An extension where all players play in a bounded convex region with a planar exit is also investigated.
	
	\item For multiple pursuers and one evader, we characterize ES by a potential function and based on ES, a guaranteed pursuer winning strategy is designed.
	\item We consider all possible coalitions among pursuers for multiple pursuers and one evader case. We further prove that in 3D if a pursuit team can guard the goal region against one evader, then at most three pursuers in the team are necessarily needed.
	\item We also consider a special subgame of degree and solve its associated HJI equation by a convex program. The HJI equation is hard to solve in general, and even harder for multiplayer cases \cite{JFF-MC-CJT-SSS:15,TM-MS-AA:17}.
	\item In order to capture the maximum number of evaders in the sense of open-loop, we propose a new class of constrained matching problems, i.e., \emph{maximum bipartite matching with conflict graph}. This problem is about assigning workers to jobs where some jobs need more than one worker to complete and one worker cannot simultaneously take more than one job. By polynomially reducing from the 3-dimensional  matching  problem, we prove that this class of constrained matching problem is NP-hard.
	\item Finally, we design the polynomial-time constant-factor Sequential Matching Algorithm to approximately solve the matching problem and show its APX-completeness. Based on the matching, a receding horizon capture strategy is proposed.
\end{enumerate}

\emph{Paper organization:} We introduce the multiplayer reach-avoid differential games in Section \ref{ProDessection}, including problem description, information structure and assumptions. Section \ref{multiplePursection} presents the main results of the case where multiple pursuers defend against one evader. In Section \ref{MaximumSec}, by solving a constrained matching problem, we design a receding horizon strategy for the pursuers to capture the maximum number of evaders in the open-loop sense. An extension to the case of a bounded convex play region with an exit is discussed in Section \ref{BoundedSec}. Numerical results are presented in Section \ref{simulationsec}, and we conclude the paper in Section \ref{conclusionsec}.

\emph{Notation:} Let $\mathbb{0}_{m\times n}$ be an $m\times n$ zero matrix. For any finite set $S$, the cardinality of $S$ is given by $|S|$, the set of non-empty subsets is given by $[S]^{+}$, and the set of non-empty subsets with cardinality less than or equal to $i$  for $1\leq i\leq|S|$ is denoted by $[S]^i$. For any subset $S$ of a topological space $X$, denote its boundary by $\partial S$. Let $\mathbb{R}$ and $\mathbb{R}^{+}$ be the set of reals and positive reals, respectively. Let $\mathbb{R}^n$ be the set of $n$-dimensional real column vectors and $\|\cdot\|_2$ be the Euclidean norm. Denote the unit sphere in $\mathbb{R}^3$ by $\mathbb{S}^2$. Let $\mathbf{x}=
[x\ \ y\ \ z]^\top\in\mathbb{R}^3$.

\section{Problem Description}\label{ProDessection}
\subsection{Multiplayer Reach-Avoid Differential Games}
\quad Consider a reach-avoid differential game with $N_p+N_e$ players, where there are $N_p$ pursuers  $\mathcal{P}=\{P_1,\dots,P_{N_p}\}$ and $N_e$ evaders $\mathcal{E}=\{E_1,\dots,E_{N_e}\}$. The players are assumed to be mass points and they have simple motion as Isaacs states \cite{RI:65}, i.e., they are holonomic. The game is played in the 3D Euclidean  space $\mathbb{R}^3$, where a plane $\mathcal{T}$ divides the game space $\mathbb{R}^3$ into two disjoint subregions $\Omega_{\rm goal}$ and $\Omega_{\rm play}$. The mathematical descriptions of $\mathcal{T}$, $\Omega_{\rm goal}$ and $\Omega_{\rm play}$ are given by $\{\mathbf{x}\in\mathbb{R}^3\,|\,z=0\},\{\mathbf{x}\in\mathbb{R}^3\,|\,z\leq0\}$ and $\{\mathbf{x}\in\mathbb{R}^3\,|\,z>0\}$, respectively. 
Let $\mathbf{x}_{P_i}(t)=[
x_{P_i}(t)\ \ y_{P_i}(t)\ \ z_{P_i}(t)]^\top\in\mathbb{R}^3$ and $\mathbf{x}_{E_j}(t)=
[
x_{E_j}(t)\ \ y_{E_j}(t) \ \ z_{E_j}(t)
]^\top\in\mathbb{R}^3
$ be the positions of $P_i$ and $E_j$ at time $t$, respectively. The dynamics of the players are described by the following decoupled systems for $t\ge0$:
\begin{equation}\label{dynamics}\begin{aligned}
\dot{\mathbf{x}}_{P_i}(t)&=v_{P_i}\mathbf{u}_{P_i}(t),&\mathbf{x}_{P_i}(0)&=\mathbf{x}^0_{P_i}, &P_i\in\mathcal{P},\\
\dot{\mathbf{x}}_{E_j}(t)&=v_{E_j}\mathbf{u}_{E_j}(t),&\mathbf{x}_{E_j}(0)&=\mathbf{x}_{E_j}^0, &E_j\in\mathcal{E},
\end{aligned}\end{equation}
where $\mathbf{x}_{P_i}^0$ and $\mathbf{x}_{E_j}^0$ are the initial positions of $P_i$ and $E_j$, and $v_{P_i}\in\mathbb{R}^+$ and $v_{E_j}\in\mathbb{R}^+$ denote the speeds of $P_i$ and $E_j$, respectively. The control inputs at time $t$ for $P_i$ and $E_j$ are their respective instantaneous headings $\mathbf{u}_{P_i}(t)$ and $\mathbf{u}_{E_j}(t)$, which satisfy the constraint $\mathbf{u}_{P_i}(t),\mathbf{u}_{E_j}(t)\in\mathbb{S}^2$. There are no other constraints on the control inputs, and all players are allowed to change their orientations instantaneously. For notational simplicity, the time $t$ will be omitted hereafter.

\quad Suppose that the pursuer $P_i$ has capture radius $r_i\ge0$. The evader $E_i$ is captured as soon as his distance from at least one of pursuers  becomes equal to the corresponding capture radius. The capture set of the pursuit team is defined by $\mathcal{C}:=\cup_{i=1}^{N_p}\mathcal{C}_i$, where $\mathcal{C}_i$ is the capture set of pursuer $P_i$ and is given by $\big\{\mathbf{x}\in\mathbb{R}^3\,|\,\|\mathbf{x}-\mathbf{x}_{P_i}\|_2\leq r_i\big\}$. Assume that the number of pursuers remains constant, and the pursuers chase the evaders until all evaders in the play region $\Omega_{\rm play}$ are captured.

\begin{figure}
	\centering
	\graphicspath{{figure_original/}}
	\subfigure{
		\includegraphics[width=80mm,height=41mm]{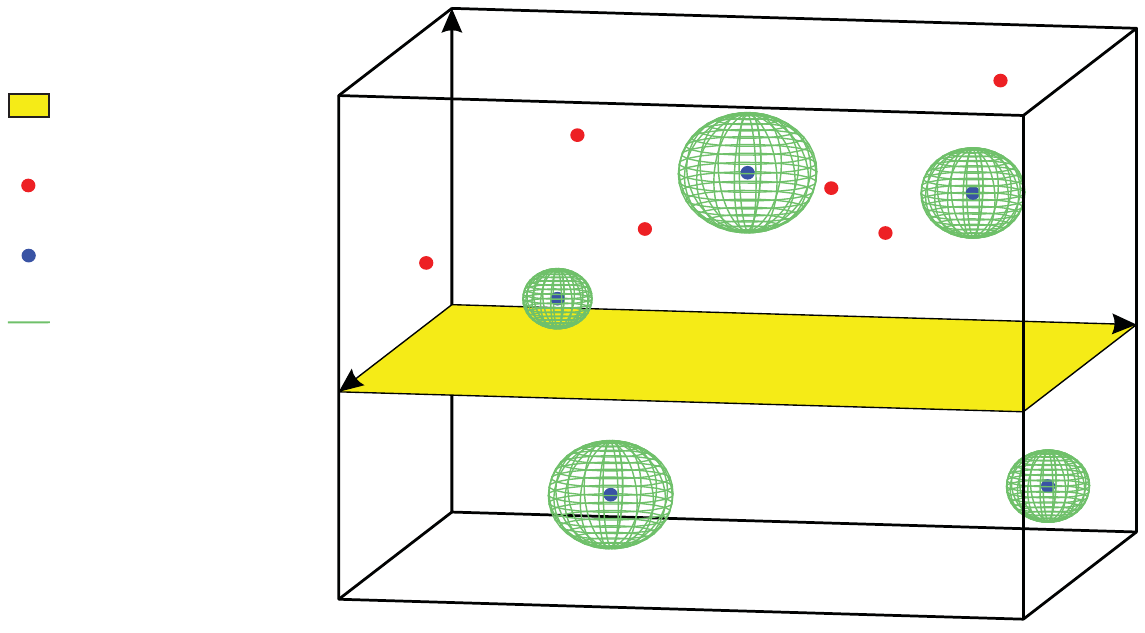}
		\put(-217,94.5){\tiny{Separating Plane}}
		\put(-217,80){\tiny{Evader}}
		\put(-217,67){\tiny{Pursuer}}
		\put(-217,55){\tiny{Capture Sphere}}
		\put(-80,10){\scriptsize{$\Omega_{\rm goal}$}}
		\put(-80,47){\scriptsize{$\mathcal{T}$}}
		\put(-80,103){\scriptsize{$\Omega_{\rm play}$}}
		\put(-120,83){\scriptsize{$E_j$}}
		\put(-88,77){\scriptsize{$P_i$}}
		\put(-72,70){\scriptsize{$r_i$}}
		\put(-157,51){\scriptsize{$x$}}
		\put(-12,60){\scriptsize{$y$}}
		\put(-135,107){\scriptsize{$z$}}
	}
	\caption{The 3D multiplayer reach-avoid differential games, where the pursuit team with multiple pursuers (blue circles) aims to capture as many evaders (red circles) as possible before these evaders penetrate a separating plane $\mathcal{T}$ (yellow) and enter the goal region $\Omega_{\rm goal}$. Each player is allowed to have different speed and each pursuer has a possibly different capture radius (green sphere). Our goal is to find a strategy for the pursuit team to capture a good number of evaders.}
	\label{game_figure}
\end{figure}
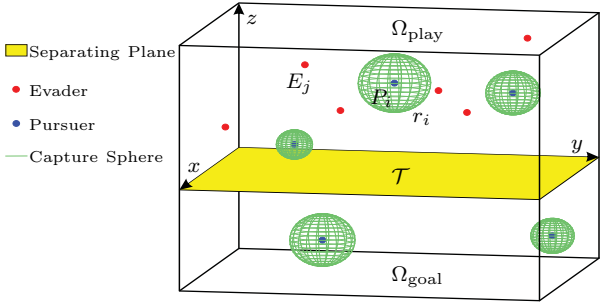

\quad The evasion team tries to send as many evaders as possible into $\Omega_{\rm goal}$ before being captured, while the pursuit team aims at capturing as many evaders as possible before they enter $\Omega_{\rm goal}$. This paper presents a receding horizon capture strategy for the pursuit team to capture as many evaders as possible. The game components are shown in Fig.~\ref{game_figure}.

\subsection{Information Structure and Assumptions}
\quad As is the usual convention in the differential game theory, the information available to each player plays an important role in generating optimal strategies \cite{TB-GJO:99}. In this paper, we adopt the state feedback information structure, where each player chooses its current input, $\mathbf{u}_{P_i}$ or $\mathbf{u}_{E_j}$, based on the current value of the information set $\{\mathbf{x}_{P_i},\mathbf{x}_{E_j}\}_{P_i\in\mathcal{P},E_j\in\mathcal{E}}$.

\quad Assume that the initial configurations of all players satisfy the following conditions, which can focus our attention on the main situations and remove some technical problems (eg., two pursuers or evaders initially lie at the same position).

\begin{asp}\label{IsolationDeploy}{\rm (Initial Deployment).}
	The initial positions of all players satisfy the following four conditions:\\
	1) $\|\mathbf{x}^0_{P_i}-\mathbf{x}^0_{P_j}\|_2>0$ for all $P_i,P_j\in\mathcal{P},P_i\neq P_j$;\\
	2) $\|\mathbf{x}^0_{E_i}-\mathbf{x}^0_{E_j}\|_2>0$ for all $E_i,E_j\in\mathcal{E},E_i\neq E_j$;\\
	3) $\|\mathbf{x}^0_{E_j}-\mathbf{x}^0_{P_i}\|_2>r_i$ for all $P_i\in\mathcal{P},E_j\in\mathcal{E}$;\\
	4) $\mathbf{x}^0_{P_i}\in\mathbb{R}^3$ for all $P_i\in\mathcal{P}$ and $\mathbf{x}^0_{E_j}\in\Omega_{\rm play}$ for all $E_j\in\mathcal{E}$.
\end{asp}
In \aspref{IsolationDeploy}, conditions $1)$ and $2)$ guarantee that all players play the game from different initial positions, condition $3)$ ensures that evaders are not captured by the pursuers initially, and condition $4)$ says that every evader initially lies in $\Omega_{\rm play}$ while every pursuer may start from any position.

\quad Most of current works on multiplayer reach-avoid differential games focus on homogeneous players in both teams \cite{PA-SDB-FB:14d,RY-ZS-YZ:19,RY-ZS-YZ:20}. We instead consider heterogeneous players, i.e., all players are allowed to have different speeds and all pursuers are allowed to have different capture radii. We focus on the faster pursuer case.
\begin{asp}\label{speedratio}{\rm (Speed Ratio).}
	Suppose the speed ratio $\alpha_{ij}=v_{P_i}/v_{E_j}>1$ for all
	$P_i\in\mathcal{P}$ and $E_j\in\mathcal{E}$.
\end{asp}

%

\subsection{A Useful Lemma}

\begin{lema}\label{polarconvexlema}{\rm (Convexity of Sets in Polar Coordinates).}
	For a twice differentiable simple 2D closed curve $\rho:[0,2\pi]\mapsto\mathbb{R}^{+}$ with $\rho(0)=\rho(2\pi)$, the set of points consisting of this curve and its interior is strictly convex if for all $\psi\in[0,2\pi]$, $\rho(\psi)$ satisfies
	\begin{equation}\label{eq:curvature}
	\rho^2+2\big(\frac{d \rho}{d \psi}\big)^2-\rho\frac{d ^2\rho}{d \psi^2}>0.
	\end{equation}
	Proof. \rm
		We postpone the proof to Appendix~\ref{appx:convexcurve}. \qed
\end{lema}

\section{Multiple Pursuers Versus One Evader}\label{multiplePursection}
\subsection{Problem statement}
\quad Since it is hard to analyze the whole game directly \cite{MC-ZZ-CJT:17}, we will decompose the whole game as many subgames which involve multiple pursuers and only one evader. Then, these subgames are used as a building block for the pursuer-evader matching in the next section. In \cite{MC-ZZ-CJT:17}, the authors only consider the subgames between one pursuer and one evader before the matching. Next, these subgames will be discussed.

\quad For any $s\in\big[\{1,\dots,N_p\}\big]^{+}$, let $P_s=\{P_i\in\mathcal{P}\,|\,i\in s\}$ be an element of $[\mathcal{P}]^{+}$, and we refer to $P_s$ as a pursuit coalition containing  pursuer $P_i$ if the subscript satisfies $i\in s$. In other word, $P_s$ is a pursuit coalition with its members specified by the index set $s$. For $P_s\in[\mathcal{P}]^{+}$, we stack the states and control inputs of all pursuers in $P_s$ into $\mathbf{x}_s$ and $\mathbf{u}_s$ respectively. Let $\mathbf{x}_s^0$ be the initial state of $\mathbf{x}_s$. 

\quad For the subgame between a pursuit coalition $P_s$ and an evader $E_j$, the pursuit coalition $P_s$ wins if it can capture $E_j$ before the latter reaches $\Omega_{\rm goal}$, while the evader $E_j$ wins if it can reach $\Omega_{\rm goal}$ before being captured by $P_s$. The capture means $\mathbf{x}_{E_j}\in\mathcal{C}_s$, where $\mathcal{C}_s:=\cup_{i\in s}\mathcal{C}_i$. 

\quad This section will address the following problems.

\begin{pbm}
	Consider a pursuit coalition $P_s$ and an evader $E_j$. Given $\mathbf{x}_s$ and $\mathbf{x}_{E_j}$, find the conditions to determine which one can win the game. If $P_s$ can win the game, what are the strategies to be adopted by the pursuers in $P_s$ to guarantee this winning? 
\end{pbm}

\begin{pbm}
	Consider a pursuit coalition $P_s$ and an evader $E_j$. If $P_s$ can win against $E_j$, how many pursuers at most in $P_s$ are necessary to guarantee this winning?   
\end{pbm}

\subsection{Evasion Space}\label{ESsection}
\quad We consider the concept introduced in Section 6.7 in Isaacs' book \cite{RI:65} as follows. 

\begin{defi}\label{ESDefi}{\rm (Evasion Space).}
	Given any $P_s\in[\mathcal{P}]^{+}$ and $E_j\in\mathcal{E}$, the evasion space (ES) $\mathbb{E}(s,j)$ is the set of positions in $\mathbb{R}^3$ that $E_j$ can reach without being captured by $P_s$, regardless of $P_s$' control input, and let the surface $\partial \mathbb{E}(s,j)$ which bounds the space $\mathbb{E}(s,j)$ be designated by the BES (boundary of evasion space).
\end{defi}

\quad First, we introduce a class of potential functions.
\begin{defi}\label{potentialfunc}{\rm (Potential Function).}
	Given $\mathbf{x}_{P_i}$ and $\mathbf{x}_{E_j}$ satisfying $\|\mathbf{x}_{E_j}-\mathbf{x}_{P_i}\|_2>r_i$, define the potential function $f_{ij}(\mathbf{x}):\mathbb{R}^3\mapsto\mathbb{R}$ associated with $P_i$ and $E_j$ as follows
	\begin{equation}\label{PFexpression}
	f_{ij}(\mathbf{x})=\|\mathbf{x}-\mathbf{x}_{P_i}\|_2-\alpha_{ij}\|\mathbf{x}-\mathbf{x}_{E_j}\|_2-r_i,
	\end{equation}
	whose gradient with respect to $\mathbf{x}$ is denoted by $\nabla f_{ij}(\mathbf{x})\in\mathbb{R}^3$, and given by
	\begin{equation}\label{FDerivative}
	\nabla f_{ij}(\mathbf{x})=\frac{\mathbf{x}-\mathbf{x}_{P_i}}{\|\mathbf{x}-\mathbf{x}_{P_i}\|_2}-\alpha_{ij}\frac{\mathbf{x}-\mathbf{x}_{E_j}}{\|\mathbf{x}-\mathbf{x}_{E_j}\|_2},
	\end{equation}
	when $\mathbf{x}\neq\mathbf{x}_{P_i}$ and $\mathbf{x}\neq\mathbf{x}_{E_j}$.
\end{defi}

\begin{figure}
	\centering
	\graphicspath{{figure_original/}}
	\subfigure{
		\includegraphics[width=85mm,height=30mm]{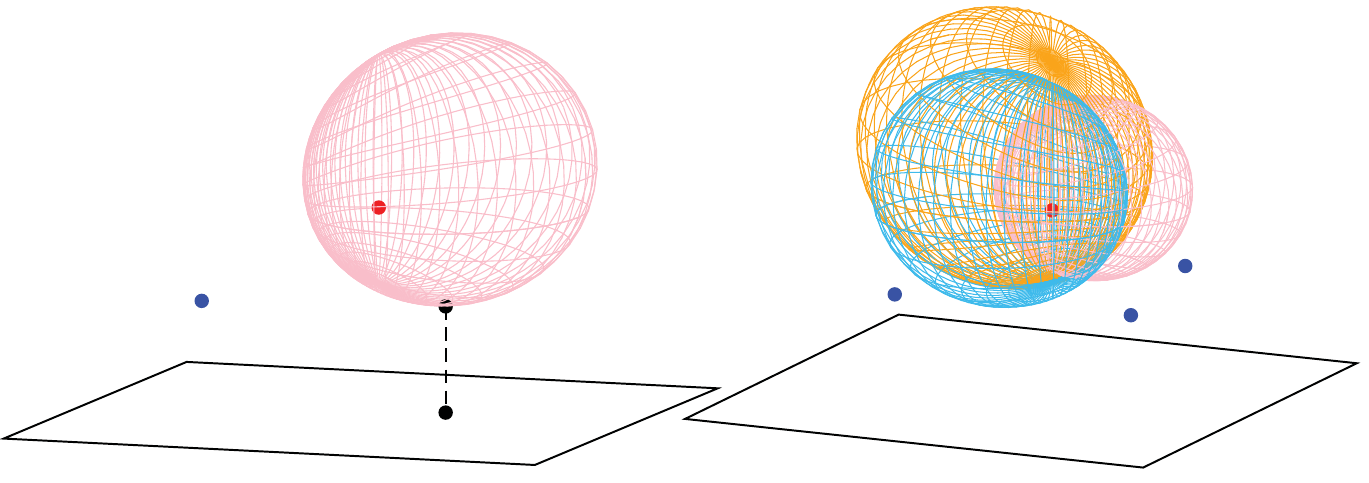}
		\put(-185,53){\scriptsize{$E_j$}}
		\put(-214,36){\scriptsize{$P_i$}}
		\put(-146,75){\scriptsize{$\bar{\mathbb{E}}(i,j)$}}
		\put(-161,25){\scriptsize{$I(i,j)$}}
		\put(-223,9){\scriptsize{$\mathcal{T}$}}
		\put(-212,-2){\scriptsize{$\Omega_{\rm goal}$}}
		\put(-173,18.8){{$\Big\{$}}
		\put(-194,23){\scriptsize{$z_{I(i,j)}$}}
		\put(-96,30.5){\scriptsize{$P_{1}$}}
		\put(-48,19.5){\scriptsize{$P_{2}$}}
		\put(-31,34){\scriptsize{$P_{3}$}}
		\put(-67,52){\scriptsize{$E_j$}}
		\put(-107,12){\scriptsize{$\mathcal{T}$}}
		\put(-98,0){\scriptsize{$\Omega_{\rm goal}$}}
		\put(-48.5,79){\scriptsize{$\bar{\mathbb{E}}(2,j)$}}
		\put(-30,50){\scriptsize{$\bar{\mathbb{E}}(1,j)$}}
		\put(-108.5,39){\scriptsize{$\bar{\mathbb{E}}(3,j)$}}
		\put(-185,-5){\scriptsize{$(a)$}}
		\put(-60,-5){\scriptsize{$(b)$}}
	}
	\caption{Evasion space (ES) and interception point. $(a)$ Single-pursuer case: The closure of ES $\bar{\mathbb{E}}(i,j)$ (pink) associated with a pursuer $P_i$ and an evader $E_j$ is a strictly convex set, and the related interception point $I(i,j)$ is the unique point in $\bar{\mathbb{E}}(i,j)$ that is closest to $\Omega_{\rm goal}$. $(b)$ Multiple-pursuer case: The closure of ES $\bar{\mathbb{E}}(s,j)$ associated with a pursuit coalition $P_s=\{P_1,P_2,P_3\}$ and an evader $E_j$ is the intersection set of three one-to-one ESs, i.e., $\cap_{i\in s}\bar{\mathbb{E}}(i,j)$. Thus, $\bar{\mathbb{E}}(s,j)$ is strictly convex and the related interception point is defined similarly which is hard to visualize here.
	}
	\label{ES_figure}
\end{figure}
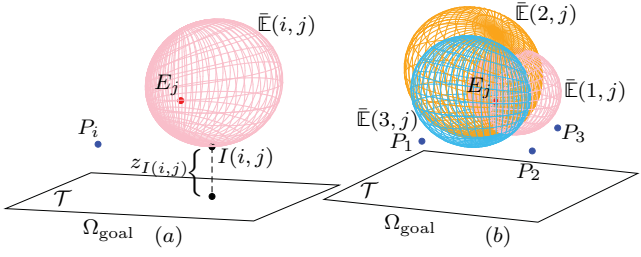

\quad Thus, according to Section 6.7 in Isaacs' book \cite{RI:65}, the ES $\mathbb{E}(i,j)$ and the BES $\partial\mathbb{E}(i,j)$ for $P_i$ and $E_j$ can be respectively given by 
\begin{equation}\label{ESBES}
\begin{aligned}
\mathbb{E}(i,j)&=\big\{\mathbf{x}\in\mathbb{R}^3\,|\,f_{ij}(\mathbf{x})>0\big\},\\
\partial\mathbb{E}(i,j)&=\big\{\mathbf{x}\in\mathbb{R}^3\,|\,f_{ij}(\mathbf{x})=0\big\}.
\end{aligned}
\end{equation}
Thus, the closure of the ES $\mathbb{E}(i,j)$, denoted by $\bar{\mathbb{E}}(i,j)$, is given by $\bar{\mathbb{E}}(i,j)=\big\{\mathbf{x}\in\mathbb{R}^3\,|\,f_{ij}(\mathbf{x})\ge0\big\}$, shown in Fig. \ref{ES_figure}(a). Next, we prove that this closure has good features.

\begin{lema}\label{ESone}{\rm (ES for One Pursuer).}
	For the ES $\mathbb{E}(i,j)$ with respect to $P_i\in\mathcal{P}$ and $E_j\in\mathcal{E}$, its closure $\bar{\mathbb{E}}(i,j)$ is bounded and strictly convex.
\end{lema}
 \emph{Proof.} \rm
	We build a new polar coordinate system with $\mathbf{x}_{E_j}$ as the origin, and let $\mathbf{x}=\mathbf{x}_{E_j}+\rho\mathbf{e}$, where $\rho\in\mathbb{R}^{+}$ and $\mathbf{e}\in\mathbb{S}^2$. We parameterize $\mathbf{e}$ by
	\begin{equation}\begin{aligned}\label{eformula}
	\mathbf{e}=\big(&\cos(\psi+\psi_0)\cos(\theta+\theta_0),\cos(\psi+\psi_0)\sin(\theta+\theta_0),\\
	&\sin(\psi+\psi_0)\big),\\
	&\theta\in[0,\pi],\theta_0\in[0,\pi],\psi\in[0,2\pi],\psi_0\in[0,2\pi],
	\end{aligned}\end{equation}
	where $\theta$ and $\psi$ are rotations with respect to positive $x$-axis and $x$-$y$-plane respectively, and $\theta_0$ and $\psi_0$ are initial rotations with respect to the original coordinates. By \eqref{PFexpression} and \eqref{ESBES}, the boundary $\partial\mathbb{E}(i,j)$, i.e., $f_{ij}(\mathbf{x})=0$, in this polar coordinates becomes
	\begin{equation}\begin{aligned}\label{polarcoor}
	&\|\mathbf{x}_{E_j}+\rho\mathbf{e}-\mathbf{x}_{P_i}\|_2-\alpha_{ij}\rho-r_i=0\\
	&\Rightarrow\rho=\frac{1}{\alpha_{ij}^2-1}\big(h_1(\theta,\psi)+h_2(\theta,\psi)\big),
	\end{aligned}\end{equation}
	where two scalar functions $h_1(\theta,\psi)$ and $h_2(\theta,\psi)$ are
	\begin{equation}\begin{aligned}\label{h1h2defi}
	h_1(\theta,\psi)&=(\mathbf{x}_{E_j}-\mathbf{x}_{P_i})^\top\mathbf{e}-\alpha_{ij}r_i,\\
	h_2(\theta,\psi)&=\sqrt{h_1^2(\theta,\psi)+(\alpha_{ij}^2-1)(\|\mathbf{x}_{E_j}-\mathbf{x}_{P_i}\|_2^2-r_i^2)}.
	\end{aligned}
	\end{equation}
	In deriving \eqref{polarcoor} and \eqref{h1h2defi}, we have used the fact that $\alpha_{ij}>1$ and $\|\mathbf{x}_{E_j}-\mathbf{x}_{P_i}\|_2>r_i$, which also implies that $h_2>0$. From \eqref{polarcoor}, we have that $\rho$ is bounded and $\rho>0$, and thus $\bar{\mathbb{E}}(i,j)$ is bounded.
	
	\quad Regarding the strict convexity of $\bar{\mathbb{E}}(i,j)$, according to the definition of $h_2(\theta,\psi)$ in \eqref{h1h2defi}, we have
	\begin{equation}\begin{aligned}\label{h1h2relation}
	\frac{\partial h_2}{\partial \psi}&=\frac{h_1}{\sqrt{h_1^2+(\alpha_{ij}^2-1)(\|\mathbf{x}_{E_j}-\mathbf{x}_{P_i}\|_2^2-r_i^2)}}\frac{\partial h_1}{\partial \psi}\\
	&=\frac{h_1}{h_2}\frac{\partial h_1}{\partial \psi}.
	\end{aligned}\end{equation}
	It follows from \eqref{polarcoor} that by fixing $\theta$, the first and second order partial derivatives of $\rho$ with respect to $\psi$ are
	\begin{equation}
	\begin{aligned}\label{twoderiveativeque}
	\frac{\partial \rho}{\partial \psi}&=\frac{1}{\alpha_{ij}^2-1}\frac{h_2+h_1}{h_2}\frac{\partial h_1}{\partial \psi},\\
	\frac{\partial^2 \rho}{\partial \psi^2}&=\frac{1}{\alpha_{ij}^2-1}\Big(\frac{h_2+h_1}{h_2}\frac{\partial^2 h_1}{\partial \psi^2}+\frac{h_2^2-h_1^2}{h_2^3}\big(\frac{\partial h_1}{\partial \psi}\big)^2\Big),
	\end{aligned}
	\end{equation}
	where \eqref{h1h2relation} is used. Then, by \eqref{polarcoor} and \eqref{twoderiveativeque}, we have
	\[\begin{aligned}
	&\rho^2+2\big(\frac{\partial \rho}{\partial \psi}\big)^2-\rho\frac{\partial ^2\rho}{\partial \psi^2}\\
	&=\frac{(h_1+h_2)^2}{(\alpha_{ij}^2-1)^2}+2\frac{1}{(\alpha_{ij}^2-1)^2}\Big(\frac{h_2+h_1}{h_2}\frac{\partial h_1}{\partial \psi}\Big)^2\\
	&\quad-\frac{h_1+h_2}{\alpha_{ij}^2-1}\frac{1}{\alpha_{ij}^2-1}\Big(\frac{h_2+h_1}{h_2}\frac{\partial^2 h_1}{\partial \psi^2}+\frac{h_2^2-h_1^2}{h_2^3}\big(\frac{\partial h_1}{\partial \psi}\big)^2\Big)\\
	&=\frac{(h_2+h_1)^2}{(\alpha_{ij}^2-1)^2}\Big(1-\frac{1}{h_2}\frac{\partial^2 h_1}{\partial \psi^2}+\frac{h_2+h_1}{h_2^3}\big(\frac{\partial h_1}{\partial \psi}\big)^2\Big)\\
	&\ge\frac{(h_2+h_1)^2}{(\alpha_{ij}^2-1)^2}\Big(1-\frac{1}{h_2}\frac{\partial^2 h_1}{\partial \psi^2}\Big).
	\end{aligned}\]
	
	\quad In the following, we prove that $h_2>\frac{\partial^2 h_1}{\partial \psi^2}$. Note that by \eqref{h1h2defi}, $h_2>0$ and $\frac{\partial^2 h_1}{\partial \psi^2}=-(\mathbf{x}_{E_j}-\mathbf{x}_{P_i})^\top\mathbf{e}=-h_1-\alpha_{ij}r_i$. On the one hand, if $h_1>-\alpha_{ij}r_i$, then we have $h_2>0>\frac{\partial^2 h_1}{\partial \psi^2}$. On the other hand, if $h_1\leq-\alpha_{ij}r_i$, then we have
	\begin{align*}
	h_2^2&-\big(\frac{\partial^2 h_1}{\partial \psi^2}\big)^2\\
	&=(\alpha_{ij}^2-1)(\|\mathbf{x}_{E_j}-\mathbf{x}_{P_i}\|_2^2-r_i^2)-\alpha_{ij}r_i(2h_1+\alpha_{ij}r_i)\\
	&>-\alpha_{ij}r_i(2h_1+\alpha_{ij}r_i)\ge \alpha_{ij}^2r_i^2\ge0.
	\end{align*}
	
	\quad Thus, we have $\rho^2+2\big(\frac{\partial \rho}{\partial \psi}\big)^2-\rho\frac{\partial ^2\rho}{\partial \psi^2}>0$ for all $\psi$. It follows from \lemaref{polarconvexlema} that $\bar{\mathbb{E}}(i,j)$ is strictly convex for any fixed $\theta$. Also note that we can take any $\theta_0\in[0,\pi]$ and $\psi_0\in[0,2\pi]$ as the initial rotations. Thus, by taking any admissible $\theta_0$ and $\psi_0$ and then considering all $\theta$ in $[0,\pi]$, we obtain that $\bar{\mathbb{E}}(i,j)$ is strictly convex. Here, it is worth emphasizing again that the strict convexity is correct because any admissible $\theta_0$ and $\psi_0$ can be taken as initial rotations. \qed

\quad It can be seen that if $r_i=0$, the BES $\partial\mathbb{E}(i,j)$ is actually the Apollonius circle in 3D (see Section 6.7 in Isaacs' book \cite{RI:65}) proposed for the case of zero capture radius. However, for the case of non-zero capture radius, the analysis is far beyond the scope of Apollonius circle. Most of current works on multiplayer reach-avoid differential games focus on zero capture radius in terms of analytical results, because the non-zero capture radius greatly increases the complexity of explicit computation \cite{PA-SDB-FB:14d,DS-VK:18,RY-ZS-YZ:20,EG-DWC-MP:19}. Recently, a variation of Apollonius circle is introduced in \cite{RY-ZS-YZ:19-3}, based on a given time difference between two players' minimal arrival time to a point. 

\quad Next we compute the ES when multiple pursuers are involved, as Fig. \ref{ES_figure}(b) shows. According to Section 6.8 in Isaacs' book \cite{RI:65} and \defiref{ESDefi}, the ES $\mathbb{E}(s,j)$ with respect to a pursuit team $P_s$ and an evader $E_j$ can be computed by
\[
\mathbb{E}(s,j)=\big\{\mathbf{x}\in\mathbb{R}^3\,|\,f_{ij}(\mathbf{x})>0,i\in s\big\}.
\]
Thus, the closure of $\mathbb{E}(s,j)$, denoted by $\bar{\mathbb{E}}(s,j)$, is given by $\bar{\mathbb{E}}(s,j)=\big\{\mathbf{x}\in\mathbb{R}^3\,|\,f_{ij}(\mathbf{x})\ge0,i\in s\big\}$. Also notice that $\mathbb{E}(s,j)=\cap_{i\in s}\mathbb{E}(i,j)$ and $\partial\mathbb{E}(s,j)\subseteq\cup_{i\in s}\partial\mathbb{E}(i,j)$. According to \lemaref{ESone} and the properties of the intersection set of convex sets,  $\bar{\mathbb{E}}(s,j)$ is also bounded and strictly convex.


\subsection{Guaranteed Pursuer Winning Strategies}
\quad This subsection considers the case when $\bar{\mathbb{E}}(s,j)\cap\Omega_{\rm goal}$ is empty, implying that there is no point in $\Omega_{\rm goal}$ that $E_j$ can reach without being captured by $P_s$. We will propose a feedback strategy for the pursuit coalition $P_s$ to guarantee its winning in this case. Note that the closure of evasion space $\bar{\mathbb{E}}(s,j)$ is bounded and strictly convex.

\begin{defi}\label{goalpointlema}{\rm (Interception Point).}
	For $P_s\in[\mathcal{P}]^{+}$ and $E_j\in\mathcal{E}$, if $\bar{\mathbb{E}}(s,j)\cap\Omega_{\rm goal}$ is empty, let the interception point $I(s,j)$ $=
	\big[x_{I(s,j)}\ y_{I(s,j)}\ z_{I(s,j)}\big]^\top
	\in\bar{\mathbb{E}}(s,j)$ be the unique point that is closest to $\Omega_{\rm goal}$.
\end{defi}

The geometrical meaning of the interception point $I(s,j)$ is shown in Fig. \ref{ES_figure}. The interception point has the following properties.
\begin{lema}\label{GPBESpro}{\rm (Properties of the Interception Point).}
	Given any $P_s\in[\mathcal{P}]^{+}$ and $E_j\in\mathcal{E}$, suppose that $\bar{\mathbb{E}}(s,j)\cap\Omega_{\rm goal}$ is empty. The interception point $I(s,j)$ has the following properties:
	\begin{enumerate}[(i)]
		\item $I(s,j)$ lies on $\partial\mathbb{E}(s,j)$;
		\item for any $s$ with $|s|=3$, if $E_j$ and the pursuers in $P_s$ are not coplanar and $I(s,j)\in\bigcap_{i\in s}\partial\mathbb{E}(i,j)$, then there exists a plane such that $I(s,j)$ is an intersection point of two strictly convex closed curves in the plane;
		\item for any $s$ with $|s|=3$, if $E_j$ and the pursuers in $P_s$ are not coplanar, then
		$\bigcap_{i\in s}\partial\mathbb{E}(i,j)$ contains at most $4$ intersection points.
	\end{enumerate}
	Proof. \rm
		Regarding \emph{(i)}, it follows from the strict convexity of $\bar{\mathbb{E}}(s,j)$ and the definition of $I(s,j)$.
		
		\quad Regarding \emph{(ii)}, let $P_s=\{P_1,P_2,P_3\}$, and then it follows from $I(s,j)\in\bigcap_{i\in s}\partial\mathbb{E}(i,j)$ and \eqref{ESBES}  that the interception point $I(s,j)=\mathbf{x}_{E_j}+\rho\mathbf{e}$ satisfies the following system of equations
		\begin{equation}
		\begin{aligned}\label{uniqueequ1}
		\|\mathbf{x}_{E_j}-\mathbf{x}_{P_i}+\rho\mathbf{e}\|_2=\alpha_{ij}\rho+r_i,\quad\forall i\in s,
		\end{aligned}
		\end{equation}
		where $\rho\in\mathbb{R}^{+}$ and $\mathbf{e}\in\mathbb{S}^2$. Equivalently, we have
		\begin{equation}
		\begin{aligned}\label{uniqueequ2}
		\rho^2=c_i+\rho(\mathbf{m}_i^\top\mathbf{e}-b_i),\quad\forall i\in s,
		\end{aligned}
		\end{equation}
		where $\mathbf{m}_i=2(\mathbf{x}_{E_j}-\mathbf{x}_{P_i})/(\alpha^2_{ij}-1),b_i=2\alpha_{ij}r_i/(\alpha^2_{ij}-1)$, and $c_i=(\|\mathbf{x}_{E_j}-\mathbf{x}_{P_i}\|_2^2-r_i^2)/(\alpha^2_{ij}-1)$. When we eliminate the term $\rho^2$ in~\eqref{uniqueequ2}, we then have
		\begin{equation*}
		\label{uniqueequ3}
		\begin{cases}
		\big((\mathbf{m}_1-\mathbf{m}_2)^\top\mathbf{e}+b_2-b_1\big)\rho=c_2-c_1\\
		\big((\mathbf{m}_2-\mathbf{m}_3)^\top\mathbf{e}+b_3-b_2\big)\rho=c_3-c_2
		\end{cases}
		\end{equation*}
		which implies that
		\begin{multline}\label{uniqueplanepro}
		\big((c_2-c_3)\mathbf{m}_1+(c_3-c_1)\mathbf{m}_2+(c_1-c_2)\mathbf{m}_3\big)^\top\mathbf{e}\\
		=(c_3-c_2)(b_2-b_1)-(c_2-c_1)(b_3-b_2).
		\end{multline}
		Since the four players are not coplanar, the vectors $\mathbf{m}_1,\mathbf{m}_2$ and $\mathbf{m}_3$ are linearly independent. Hence, by \eqref{uniqueplanepro}, the vector $\mathbf{e}$ lies in a plane, and thus the same for $I(s,j)$. To solve \eqref{uniqueequ1}, we could replace the case of $i=3$ in \eqref{uniqueequ1} with \eqref{uniqueplanepro}. Note that the intersection of (\ref{uniqueplanepro}) and $\partial\mathbb{E}(1,j)$ is a strictly convex closed curve, and the same for that of  (\ref{uniqueplanepro}) and $\partial\mathbb{E}(2,j)$. Thus, $I(s,j)$ is an intersection point of two strictly convex closed curves in the plane given by (\ref{uniqueplanepro}).
		
		\quad Regarding \emph{(iii)}, we show that there are at most four solutions to  \eqref{uniqueequ1}, i.e., \eqref{uniqueequ2}. We rewrite (\ref{uniqueequ2}) as $\mathbf{m}_i^\top\mathbf{e}=(\rho^2-c_i)/\rho+b_i$ for $i=1,2,3$. Then, given $\rho$, the vector $\mathbf{e}$ is uniquely given by
		\begin{equation}\label{eexpression}
		\mathbf{e}=
		\begin{bmatrix}
		\mathbf{m}_1^\top\\
		\mathbf{m}^\top_2 \\ \mathbf{m}^\top_3
		\end{bmatrix}^{-1}
		\begin{bmatrix}
		(\rho^2-c_1)/\rho+b_1\\(\rho^2-c_2)/\rho+b_2\\(\rho^2-c_3)/\rho+b_3
		\end{bmatrix},
		\end{equation}
		where the matrix inversion is well-defined because $\mathbf{m}_1$, $\mathbf{m_2}$ and $\mathbf{m}_3$ are linearly independent as the four players are not coplanar. By the fact that $\|\mathbf{e}\|_2=1$,  \eqref{eexpression} becomes a quartic equation of $\rho$, which has at most four solutions. \qed
\end{lema}

\begin{figure}
	\centering
	\graphicspath{{figure_original/}}
	\subfigure{
		\includegraphics[width=82mm,height=32mm]{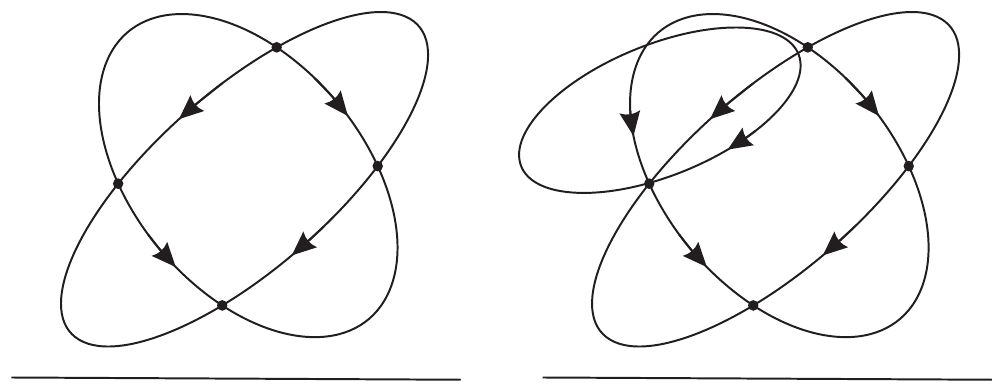}
		\put(-230,3){\scriptsize{$\mathcal{T}$}}
		\put(-105,3){\scriptsize{$\mathcal{T}$}}
		\put(-150,6){\scriptsize{$\bar{\mathbb{E}}(1,j)$}}
		\put(-25,6){\scriptsize{$\bar{\mathbb{E}}(1,j)$}}
		\put(-237,35){\scriptsize{$\bar{\mathbb{E}}(2,j)$}}
		\put(-117.5,23){\scriptsize{$\bar{\mathbb{E}}(2,j)$}}
		\put(-115,78){\scriptsize{$\bar{\mathbb{E}}(4,j)$}}
		\put(-174,84){\scriptsize{$H_1$}}
		\put(-49.5,84){\scriptsize{$H_1$}}
		\put(-219,46){\scriptsize{$H_2$}}
		\put(-78,43){\scriptsize{$H_2$}}
		\put(-187,10){\scriptsize{$H_3$}}
		\put(-62,10){\scriptsize{$H_3$}}
		\put(-143,48.5){\scriptsize{$H_4$}}
		\put(-18,48.5){\scriptsize{$H_4$}}
		\put(-187,-6){\scriptsize{$(a)$}}
		\put(-62,-6){\scriptsize{$(b)$}}
	}
	\caption{Degeneration of the interception point, which shows that at most three pursuers are needed to compute the interception point for a pursuit coalition. $(a)$ Consider $P_{s_1}=\{P_1,P_2,P_3\}$ and $E_j$, and four players are not coplanar. If $I(s_1,j)\in\cap_{i\in s_1}\partial\mathbb{E}(i,j)$, then there exists a plane intersecting with two closures of ESs $\bar{\mathbb{E}}(1,j)$ and $\bar{\mathbb{E}}(2,j)$ by two strictly convex curves respectively as depicted, which are proved to allow at most four intersection points $H_1,H_2,H_3$ and $H_4$. The interception point occurs at one of these four interception points and the $z$-coordinate of the points on these two curves decreases along the arrow. Thus, $H_3$ is the interception point. $(b)$ When adding another pursuer $P_4$, an associated strictly convex curve in this plane can be obtained as depicted. If the interception point also depends on $P_4$, it is proved that at least one pursuer in $P_{s_1}$ is redundant for computing the interception point.
	}
	\label{proof3point_figure}
\end{figure}
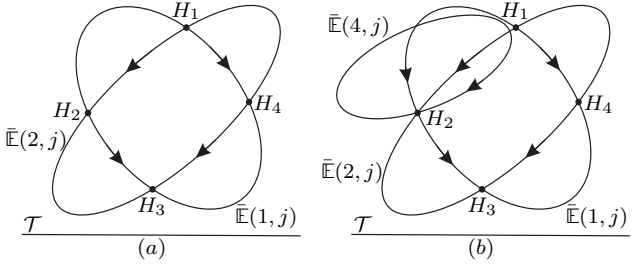

\quad Before discussing the strategies, the following
lemma is presented, which will dramatically decrease the complexity of the strategy selection and the matching in the next section.

\begin{lema}\label{degenerationlema}{\rm (Degeneration of the Interception Point).}
	For any $P_s\in[\mathcal{P}]^{+}$ and $E_j\in\mathcal{E}$, suppose that $\bar{\mathbb{E}}(s,j)\cap\Omega_{\rm goal}$ is empty. Then, there must exist a pursuit subcoalition $s_1\subseteq s$ such that $|s_1|\leq3$ and $I(s_1,j)=I(s,j)$.
	
	Proof. \rm
		The statement holds trivially if $|s|\leq3$. Therefore, we focus on the case when $|s|\ge4$. Consider $|s_1|=3$ and assume that $I(s_1,j)$ depends on all pursuers in $P_{s_1}$. Thus, $I(s_1,j)\in\cap_{i\in s_1}\partial\mathbb{E}(i,j)$. The case when $E_j$ and the three pursuers in $P_{s_1}$ are coplanar will be discussed separately below. Let $P_{s_1}=\{P_1,P_2,P_3\}$.
		
		\quad Case 1: $E_j,P_1,P_2$ and $P_3$ are not coplanar. It follows from the properties \emph{(ii)} and \emph{(iii)} in \lemaref{GPBESpro} that $I(s_1,j)$ is one of the intersection points of two strictly convex closed curves in a plane, which have at most four intersection points, as Fig. \ref{proof3point_figure}(a) shows. Note that we also replace the condition $I(s_1,j)\in\partial\mathbb{E}(3,j)$ by the fact that $I(s_1,j)$ lies in an associated plane, as the proof of the property \emph{(ii)} in \lemaref{GPBESpro} shows. The $z$-coordinate of the point on these two curves decreases along the arrow. Thus, $H_3$ is the interception point.
		
		\quad Then, we add the pursuer $P_4$ which corresponds to another strictly convex closed curve in the same plane, as Fig. \ref{proof3point_figure}(b) shows. If $I(s,j)$ depends on all four pursuers, then it must be one of the four points $H_1$, $H_2$, $H_3$ and $H_4$. If $I(s,j)=H_3$, then $P_4$ is redundant. If $I(s,j)=H_2$, as Fig. \ref{proof3point_figure}(b) indicates, then the arrow of the new curve must decrease from inside to outside of $\cap_{i\in s_1}\bar{\mathbb{E}}(i,j)$ at $H_2$; conversely, if the new curve increases from inside to outside, then $H_2$ cannot be the interception point as the point on the new curve can continue to decrease along the new curve after $H_2$ and also lies in the closure of the ES $\bar{\mathbb{E}}(s_1\cup 4,j)$. Thus, $P_2$ is redundant, because we can still conclude that $H_2$ is the interception point only by $P_1,P_3$ and $P_4$. Similarly, if $I(s,j)=H_4$, then $P_1$ is redundant. Since two curves go down at $H_1$, $I(s,j)$ cannot be $H_1$. Thus, adding a new pursuer does not increase the number of pursuers which the interception point necessarily depends on. Therefore, at most three pursuers are needed to locate the interception point.
		
		\quad Case 2: $E_j,P_1,P_2$ and $P_3$ are coplanar. Thus, the vectors $\mathbf{m}_1$, $\mathbf{m}_2$ and $\mathbf{m}_3$ are linearly dependent. Then, by following the argument of the property \emph{(ii)} in \lemaref{GPBESpro}, we can obtain (\ref{uniqueequ2}) and write  (\ref{uniqueequ2}) in the matrix form
		\begin{equation}\label{case2equ}
		\begin{bmatrix}
		\mathbf{m}_1^\top\\
		\mathbf{m}_2^\top \\
		\mathbf{m}_3^\top
		\end{bmatrix}\mathbf{e}=
		\begin{bmatrix}
		(\rho^2-c_1)/\rho+b_1\\(\rho^2-c_2)/\rho+b_2\\(\rho^2-c_3)/\rho+b_3
		\end{bmatrix},
		\end{equation}
		where $[\mathbf{m}_1 \ \ \mathbf{m}_2\ \  \mathbf{m}_3]^\top$ is singular. If \eqref{case2equ} admits solutions, then we can obtain that the first two equalities can induce the third equality directly, that is, all intersection points between $\partial\mathbb{E}(1,j)$ and $\partial\mathbb{E}(2,j)$ must belong to $\partial\mathbb{E}(3,j)$. Thus, we can ignore $P_3$ and continue to consider the remaining pursuers in $P_s$. If \eqref{case2equ} admits no solution, then there exists pursuer $P_i$ in $P_{s_1}$ such that $I(s_1,j)\notin\partial\mathbb{E}(i,j)$. Thus, we can ignore $P_i$ and continue to consider the remaining pursuers in $P_s$. \qed
\end{lema}

\quad We next present ES-based strategies. The ES-based strategies are feedback strategies, while the strategies in \cite{MC-ZZ-CJT:17} are semi-open-loop strategies.
\begin{thom}\label{ESstramultilema}{\rm (ES-Based Strategies for Multiple Pursuers).}
	For any $P_s\in[\mathcal{P}]^{+}$ and $E_j\in\mathcal{E}$, suppose that $\bar{\mathbb{E}}(s,j)\cap\Omega_{\rm goal}$ is empty, and let $s_1$ be a subset of $s$ such that $I(s_1,j)=I(s,j)$ and $|s_1|\leq3$. If every pursuer $P_i$ in $P_{s_1}$ adopts the feedback strategy $\mathbf{u}_{P_i}=\frac{I(s_1,j)-\mathbf{x}_{P_i}}{\|I(s_1,j)-\mathbf{x}_{P_i}\|_2}$, then the pursuit subcoalition $P_{s_1}$ guarantees that $\bar{\mathbb{E}}(s,j)$ does not approach $\Omega_{\rm goal}$, i.e., $\dot{z}_{I(s,j)}\ge 0$ for any $\mathbf{u}_{E_j}\in\mathbb{S}^2$. Moreover, $\dot{z}_{I(s,j)}=0$ if and only if $E_j$ adopts the feedback strategy $\mathbf{u}_{E_j}=\frac{I(s_1,j)-\mathbf{x}_{E_j}}{\|I(s_1,j)-\mathbf{x}_{E_j}\|_2}$.
\end{thom}
 \emph{Proof.} \rm
	Note that \lemaref{degenerationlema} guarantees the existence of a subcoalition $s_1$ satisfying $I(s_1,j)=I(s,j)$ and $|s_1|\leq3$. For notational convenience, let $\mathbf{x}_{I}=
	[x_I\ \ y_I\ \ z_I]^\top
	$ be the coordinate of $I(s_1,j)$. According to \defiref{goalpointlema} and the expression of $\bar{\mathbb{E}}(s_1,j)$ given in Section \ref{ESsection}, $\mathbf{x}_I$ is the unique solution of the following convex problem
	\begin{equation*}
	\begin{aligned}
	& \underset{\mathbf{x}\in\mathbb{R}^3}{\textup{minimize}}
	& & z \\
	& \textup{subject to}
	& & f_{ij}(\mathbf{x})\ge0,\quad \forall i\in s_1.
	\end{aligned}
	\end{equation*}
	The interception point $\mathbf{x}_I$ should satisfy the Karush-Kuhn-Tucker (KKT) conditions as follows:
	\begin{equation}\begin{aligned}\label{KKT}
	&[
	0\ \ 0 \ \ -1
		]^\top
	=\sum_{i\in s_1}\lambda_i\nabla f_{ij}(\mathbf{x}_I),\\
	&f_{ij}(\mathbf{x}_I)\ge0,\,\lambda_i\leq0,\,\lambda_if_{ij}(\mathbf{x}_I)=0,\quad\forall i\in s_1,
	\end{aligned}\end{equation}
	where $\lambda_i\in\mathbb{R}$ is the Lagrange multiplier. The slack conditions in \eqref{KKT} also imply that $s_1$ can be classified into two disjoint index sets $s_1^{=0}$ and $s_1^{>0}$ ($s_1^{>0}$ may be empty) satisfying
	\begin{equation}\begin{aligned}\label{slackcondition}
	\begin{cases}
	f_{ij}(\mathbf{x}_I)=0,\,\lambda_i\leq0,\quad \textup{ if }i\in s_1^{=0},\\
	f_{ij}(\mathbf{x}_I)>0,\,\lambda_i=0,\quad\textup{ if }i\in s_1^{>0}.
	\end{cases}
	\end{aligned}
	\end{equation}
	
	\quad For any $i\in s_1^{=0}$, according to \eqref{dynamics}, \eqref{PFexpression} and \eqref{FDerivative}, taking derivative of $f_{ij}(\mathbf{x}_I)=0$ at both sides with respect to $t$, we have
	\begin{equation}\begin{aligned}\label{derivativetimeequ}
	\frac{df_{ij}(\mathbf{x}_I)}{dt}=0\Rightarrow&\frac{(\mathbf{x}_I-\mathbf{x}_{P_i})^\top (\dot{\mathbf{x}}_I-v_{P_i}\mathbf{u}_{P_i})}{\|\mathbf{x}_I-\mathbf{x}_{P_i}\|_2}\\
	&=\frac{\alpha_{ij}(\mathbf{x}_I-\mathbf{x}_{E_j})^\top(\dot{\mathbf{x}}_I-v_{E_j}\mathbf{u}_{E_j})}{\|\mathbf{x}_I-\mathbf{x}_{E_j}\|_2},
	\end{aligned}\end{equation}
	namely,
	\begin{equation}\label{derivative2}
	\begin{aligned}
	&\Big(\frac{\mathbf{x}_I-\mathbf{x}_{P_i}}{\|\mathbf{x}_I-\mathbf{x}_{P_i}\|_2}-\frac{\alpha_{ij}(\mathbf{x}_I-\mathbf{x}_{E_j})}{\|\mathbf{x}_I-\mathbf{x}_{E_j}\|_2}\Big)^\top\dot{\mathbf{x}}_I\\
	&=\frac{v_{P_i}(\mathbf{x}_I-\mathbf{x}_{P_i})^\top\mathbf{u}_{P_i}}{\|\mathbf{x}_I-\mathbf{x}_{P_i}\|_2}-\frac{\alpha_{ij}v_{E_j}(\mathbf{x}_I-\mathbf{x}_{E_j})^\top\mathbf{u}_{E_j}}{\|\mathbf{x}_I-\mathbf{x}_{E_j}\|_2}.
	\end{aligned}
	\end{equation}
	We emphasize that in deriving \eqref{derivativetimeequ}, $\mathbf{x}_{P_i}$ and $\mathbf{x}_{E_j}$ in $f_{ij}(\mathbf{x}_I)$ are also functions of time $t$ as in \eqref{dynamics}.
	
	\quad By \eqref{FDerivative} and \eqref{derivative2}, the sign of the following satisfies
	\begin{equation}\begin{aligned}\label{oneESbasedinequ}
	&\nabla f_{ij}(\mathbf{x}_I)^\top\dot{\mathbf{x}}_I\\
	&=\Big(\frac{\mathbf{x}_I-\mathbf{x}_{P_i}}{\|\mathbf{x}_I-\mathbf{x}_{P_i}\|_2}-\frac{\alpha_{ij}(\mathbf{x}_I-\mathbf{x}_{E_j})}{\|\mathbf{x}_I-\mathbf{x}_{E_j}\|_2}\Big)^\top\dot{\mathbf{x}}_I\\
	&=\frac{v_{P_i}(\mathbf{x}_I-\mathbf{x}_{P_i})^\top\mathbf{u}_{P_i}}{\|\mathbf{x}_I-\mathbf{x}_{P_i}\|_2}-\frac{\alpha_{ij}v_{E_j}(\mathbf{x}_I-\mathbf{x}_{E_j})^\top\mathbf{u}_{E_j}}{\|\mathbf{x}_I-\mathbf{x}_{E_j}\|_2}\\
	&=v_{P_i}-\frac{v_{P_i}(\mathbf{x}_I-\mathbf{x}_{E_j})^\top\mathbf{u}_{E_j}}{\|\mathbf{x}_I-\mathbf{x}_{E_j}\|_2}\ge v_{P_i}-v_{P_i}=0,
	\end{aligned}\end{equation}
	where $P_i$ adopts the feedback strategy $\mathbf{u}_{P_i}=\frac{\mathbf{x}_I-\mathbf{x}_{P_i}}{\|\mathbf{x}_I-\mathbf{x}_{P_i}\|_2}$ in the second equation. Also note that \eqref{oneESbasedinequ} holds for any $\mathbf{u}_{E_j}\in\mathbb{S}^2$. Moreover, the inequality in  \eqref{oneESbasedinequ} becomes an equality if and only if $\mathbf{u}_{E_j}=\frac{\mathbf{x}_I-\mathbf{x}_{E_j}}{\|\mathbf{x}_I-\mathbf{x}_{E_j}\|_2}$. 
	
	\quad Thus, it follows from \eqref{KKT}, \eqref{slackcondition} and \eqref{oneESbasedinequ} that
	\[\begin{aligned}
	-\dot{z}_I&=
	[
	0\ \ 0 \ \ -1
	]
	\dot{\mathbf{x}}_I=\sum_{i\in s_1}\lambda_i\nabla f_{ij}(\mathbf{x}_I)^\top\dot{\mathbf{x}}_I\\
	&=\sum_{i\in s_1^{=0}}\lambda_i\nabla f_{ij}(\mathbf{x}_I)^\top\dot{\mathbf{x}}_I+\sum_{i\in s_1^{>0}}\lambda_i\nabla f_{ij}(\mathbf{x}_I)^\top\dot{\mathbf{x}}_I\\
	&=\sum_{i\in s_1^{=0}}\lambda_i\big(\nabla f_{ij}(\mathbf{x}_I)^\top\dot{\mathbf{x}}_I\big)\leq0,
	\end{aligned}\]
	holds for any $\mathbf{u}_{E_j}\in\mathbb{S}^2$ when $P_i$ adopts the strategy $\mathbf{u}_{P_i}=\frac{\mathbf{x}_I-\mathbf{x}_{P_i}}{\|\mathbf{x}_I-\mathbf{x}_{P_i}\|_2}$. Moreover, $\dot{z}_I=0$ if and only if $E_j$ adopts the feedback strategy $\mathbf{u}_{E_j}=\frac{\mathbf{x}_I-\mathbf{x}_{E_j}}{\|\mathbf{x}_I-\mathbf{x}_{E_j}\|_2}$. \qed
	
\begin{cor}{\rm  (Game of kind)}\label{Gameofkind}
	The game winner between a pursuit coalition $P_s$ and an evader $E_j$ can be determined as follows: If $\bar{\mathbb{E}}(s,j)\cap\Omega_{\rm goal}$ is empty, then the pursuit team $P_s$ wins; if $\bar{\mathbb{E}}(s,j)\cap\Omega_{\rm goal}$ has more than one element, then $E_j$ wins; if $\bar{\mathbb{E}}(s,j)\cap\Omega_{\rm goal}$ has a unique element, then two teams are tied. The maximum number of pursuers required to capture an evader before the evader reaches the goal region, is three.
	
	Proof. {\rm
		This corollary directly follows from \defiref{ESDefi} and \thomref{ESstramultilema}. \qed
	}
\end{cor}

%

\subsection{The Hamilton-Jacobi-Isaacs Equation}
\quad This subsection revisits the case when $\bar{\mathbb{E}}(s,j)\cap\Omega_{\rm goal}=\emptyset$, i.e., the pursuit coalition $P_s$ wins the game. We consider such a special subgame of degree \cite{RI:65}: Although the pursuit coalition $P_s$ can capture the evader $E_j$, the evader tries to be captured at the closest point to the goal region and the pursuit coalition seeks the opposite. Formally, the terminal set $\Psi$ and payoff function $J$ respectively are
\begin{equation}\begin{aligned}\label{gameofdegreepayoff}
&\Psi:=\big\{(\mathbf{x}_s,\mathbf{x}_{E_j})\,|\,\exists i\in s,\textup{ s.t. }\mathbf{x}_{E_j}\in\mathcal{C}_i\big\},\\
&J(\mathbf{u}_s,\mathbf{u}_{E_j};\mathbf{x}_s^0,\mathbf{x}_{E_j}^0):=z_{E_j}(t_f),
\end{aligned}
\end{equation}
where the terminal time $t_f$ is defined as the time instant when the system state enters $\Psi$. The goals of $P_s$ and $E_j$ lead to the following value function
\begin{equation}\label{gamedegreevalue}
V(\mathbf{x}_s^0,\mathbf{x}_{E_j}^0):=\max_{\mathbf{u_s}}\min_{\mathbf{u}_{E_j}}J(\mathbf{u}_s,\mathbf{u}_{E_j};\mathbf{x}_s^0,\mathbf{x}_{E_j}^0).
\end{equation}

\quad The following theorem shows that the value function $V$ for this special subgame, can be computed via a convex program for the states such that $V$ is differentiable, because in this case $V$ is the solution of an associated HJI equation \cite{RI:65}.
\begin{thom}\label{valuefunctionthom}{\rm (Value Function).}
	Consider the differential game \eqref{dynamics}, \eqref{gameofdegreepayoff} and \eqref{gamedegreevalue} where $\bar{\mathbb{E}}(s,j)\cap\Omega_{\rm goal}$ is empty. For the states $(\mathbf{x}_s,\mathbf{x}_{E_j})$ such that the value function $V(\mathbf{x}_s,\mathbf{x}_{E_j})$ is differentiable, then $V(\mathbf{x}_s,\mathbf{x}_{E_j})$ can be computed by the convex optimization problem
	\begin{equation}\label{convexpbm}
	\begin{aligned}
	V(\mathbf{x}_s,\mathbf{x}_{E_j})=\,\,& \underset{\mathbf{x}\in\mathbb{R}^3}{\textup{minimize}}
	& & z \\
	& \textup{subject to}
	& & f_{ij}(\mathbf{x})\ge0,\quad \forall i\in s.
	\end{aligned}
	\end{equation}
\end{thom}
 \emph{Proof.} \rm
	Since we only consider the states $(\mathbf{x}_s,\mathbf{x}_{E_j})$ such that the value function $V(\mathbf{x}_s,\mathbf{x}_{E_j})$ is differentiable, the value function satisfies the HJI equation \cite[Chapter 4]{RI:65} with respect to this subgame as follows:
	\begin{equation}\label{HJIequ}
	\begin{aligned}
	-\frac{\partial V(\mathbf{x}_s,\mathbf{x}_{E_j})}{\partial t}=\max_{\mathbf{u}_s}\min_{\mathbf{u}_{E_j}}&\Bigg\{\sum_{i\in s}\frac{\partial V(\mathbf{x}_s,\mathbf{x}_{E_j})^\top}{\partial\mathbf{x}_{P_i}}v_{P_i}\mathbf{u}_{P_i}\\
	&+\frac{\partial V(\mathbf{x}_s,\mathbf{x}_{E_j})^\top}{\partial\mathbf{x}_{E_j}}v_{E_j}\mathbf{u}_{E_j}\Bigg\},
	\end{aligned}
	\end{equation}
	where \eqref{dynamics} is employed. Note that in our problem $\frac{\partial V}{\partial t}=0$. Therefore, \eqref{HJIequ} can also be equivalently rewritten as
	\begin{equation}\label{newvalue}
	0=\max_{\mathbf{u}_s}\min_{\mathbf{u}_{E_j}}\frac{dV(\mathbf{x}_s,\mathbf{x}_{E_j})}{dt}=\max_{\mathbf{u}_s}\min_{\mathbf{u}_{E_j}}\dot{V}(\mathbf{x}_s,\mathbf{x}_{E_j}).
	\end{equation}
	
	\quad Next, we show that the unique optimal value to the convex optimization problem \eqref{convexpbm} satisfies \eqref{newvalue}. Let $\mathbf{x}_I$ be the solution of \eqref{convexpbm} and suppose that $V(\mathbf{x}_s,\mathbf{x}_{E_j})=z_I$. Since $\mathbf{x}_I$ satisfies the KKT conditions of \eqref{convexpbm}\begin{equation}\begin{aligned}\label{KKT3} 
	&[
	0\ \ 0\ \ -1
	]^\top
	=\sum_{i\in s}\lambda_i\nabla f_{ij}(\mathbf{x}_I),\\
	&f_{ij}(\mathbf{x}_I)\ge0,\,\lambda_i\leq0,\,\lambda_if_{ij}(\mathbf{x}_I)=0,\quad\forall i\in s,
	\end{aligned}\end{equation}
	we have
	\begin{equation}\begin{aligned}\label{compute1}
	\dot{V}(\mathbf{x}_s,\mathbf{x}_{E_j})=\dot{z}_I=[
	0\ \ 0\ \ 1
	]
	\dot{\mathbf{x}}_I=-\sum_{i\in s}\lambda_i\nabla f_{ij}(\mathbf{x}_I)^\top\dot{\mathbf{x}}_I.
	\end{aligned}\end{equation}
	The slackness conditions in \eqref{KKT3} imply that $s$ can be classified into two disjoint index sets $s^{=0}$ and $s^{>0}$ satisfying
	\begin{equation}\begin{aligned}\label{slackcondition2}
	\begin{cases}
	f_{ij}(\mathbf{x}_I)=0,\,\lambda_i\leq0,\quad \textup{ if }i\in s^{=0},\\
	f_{ij}(\mathbf{x}_I)>0,\,\lambda_i=0,\quad\textup{ if }i\in s^{>0}.
	\end{cases}
	\end{aligned}
	\end{equation} 
	
	\quad Moreover, for any $i\in s^{=0}$, similar to the argument \eqref{derivativetimeequ} and  \eqref{derivative2}, we can obtain
	\begin{equation}\label{newderivative2}
	\begin{aligned}
	&\nabla f_{ij}(\mathbf{x}_I)^\top\dot{\mathbf{x}}_I\\
	&=\frac{v_{P_i}(\mathbf{x}_I-\mathbf{x}_{P_i})^\top\mathbf{u}_{P_i}}{\|\mathbf{x}_I-\mathbf{x}_{P_i}\|_2}-\frac{\alpha_{ij}v_{E_j}(\mathbf{x}_I-\mathbf{x}_{E_j})^\top\mathbf{u}_{E_j}}{\|\mathbf{x}_I-\mathbf{x}_{E_j}\|_2}.
	\end{aligned}
	\end{equation}
	where \eqref{FDerivative} is used in the left side of \eqref{newderivative2}.
	By \eqref{compute1}, \eqref{slackcondition2}  and \eqref{newderivative2}, we compute
	\begin{equation*}\begin{aligned}
	&\max_{\mathbf{u}_s}\min_{\mathbf{u}_{E_j}}\dot{V}(\mathbf{x}_s,\mathbf{x}_{E_j})
	=\max_{\mathbf{u}_s}\min_{\mathbf{u}_{E_j}}-\sum_{i\in s}\lambda_i\nabla f_{ij}(\mathbf{x}_I)^\top\dot{\mathbf{x}}_I\\
	&=\max_{\mathbf{u}_s}\min_{\mathbf{u}_{E_j}}-\sum_{i\in s^{=0}}\lambda_i\nabla f_{ij}(\mathbf{x}_I)^\top\dot{\mathbf{x}}_I\\
	&=\max_{\mathbf{u}_s}\min_{\mathbf{u}_{E_j}}\sum_{i\in s^{=0}}\lambda_i\Big(\frac{-v_{P_i}(\mathbf{x}_I-\mathbf{x}_{P_i})^\top\mathbf{u}_{P_i}}{\|\mathbf{x}_I-\mathbf{x}_{P_i}\|_2}\\
	&\qquad\qquad\qquad\qquad\ \ +\frac{\alpha_{ij}v_{E_j}(\mathbf{x}_I-\mathbf{x}_{E_j})^\top\mathbf{u}_{E_j}}{\|\mathbf{x}_I-\mathbf{x}_{E_j}\|_2}\Big)\\
	&=\max_{\mathbf{u}_s}\sum_{i\in s^{=0}}\lambda_i\frac{-v_{P_i}(\mathbf{x}_I-\mathbf{x}_{P_i})^\top\mathbf{u}_{P_i}}{\|\mathbf{x}_I-\mathbf{x}_{P_i}\|_2}\\
	&\quad\,+\min_{\mathbf{u}_{E_j}}\sum_{i\in s^{=0}}\lambda_i\frac{\alpha_{ij}v_{E_j}(\mathbf{x}_I-\mathbf{x}_{E_j})^\top\mathbf{u}_{E_j}}{\|\mathbf{x}_I-\mathbf{x}_{E_j}\|_2}\\
	&=-\sum_{i\in s^{=0}}\lambda_iv_{P_i}+\sum_{i\in s^{=0}}\lambda_i\alpha_{ij}v_{E_j}=0,
	\end{aligned}
	\end{equation*}
	where in the max and min operations, we take $\mathbf{u}_{P_i}=\frac{\mathbf{x}_I-\mathbf{x}_{P_i}}{\|\mathbf{x}_I-\mathbf{x}_{P_i}\|_2}$ for all $i\in s^{=0}$ and $\mathbf{u}_{E_j}=\frac{\mathbf{x}_I-\mathbf{x}_{E_j}}{\|\mathbf{x}_I-\mathbf{x}_{E_j}\|_2}$ by noting that $\lambda_i\leq0$. Thus, the value function $V(\mathbf{x}_s,\mathbf{x}_{E_j})=z_I$ satisfies the HJI equation \eqref{newvalue}, i.e., \eqref{HJIequ}.
	
	\quad Finally, we prove that the terminal condition $V(\mathbf{x}_s,\mathbf{x}_{E_j})=z_{E_j}$ from \eqref{gameofdegreepayoff} is satisfied. By the definition of $\Psi$, when the game ends, there exists one pursuer $P_i$ in $P_s$ such that $\|\mathbf{x}_{E_j}-\mathbf{x}_{P_i}\|_2=r_{i}$, which implies that $\bar{\mathbb{E}}(i,j)$ contains a unique point $\mathbf{x}_{E_j}$. In other words, we have that $f_{ij}(\mathbf{x})\ge0$ leads to $\mathbf{x}=\mathbf{x}_{E_j}$. Note that for the other pursuers in $P_s$, the constraints in \eqref{convexpbm} are feasible. Thus, the convex optimization problem \eqref{convexpbm} admits the unique solution $\mathbf{x}_I=\mathbf{x}_{E_j}$. Therefore, the value function satisfies $V(\mathbf{x}_s,\mathbf{x}_{E_j})=z_I=z_{E_j}$. \qed

\begin{rek}
	\thomref{valuefunctionthom} shows that the HJI equation of a special subgame between a pursuit coalition and an evader, which describes the value function and often is hard to solve, can be transformed into a convex optimization problem with greatly reduced  computational complexity. For the strategies of the players, they are the gradients of the value function with respect to states, and in our case they can be obtained through the optimal solution to a convex optimization problem, as the proof of \thomref{valuefunctionthom} indicates.
\end{rek}

\section{Maximum-Matching Capture Strategies}\label{MaximumSec}
\subsection{Maximum Matching}
\quad We piece together the outcomes of all pursuit coalitions and evader pairs using maximum matching as in \cite{MC-ZZ-CJT:17,RY-ZS-YZ:20}. Interestingly, thanks to \corref{Gameofkind}, the matching problem is  simplified greatly as we only need to consider all pursuit coalitions of size less than or equal to three. The pursuit team $\mathcal{P}$ consists of $N_p(N_p^2+5)/6$ possible coalitions: $N_p$ one-pursuer coalitions, $N_p(N_p-1)/2$ two-pursuer coalitions, and $N_p(N_p-1)(N_p-2)/6$ three-pursuer coalitions. For notational convenience, we define the number of possible vertices for $\mathcal{P}$ in the bipartite graph by $N_o=N_p(N_p^2+5)/6$.

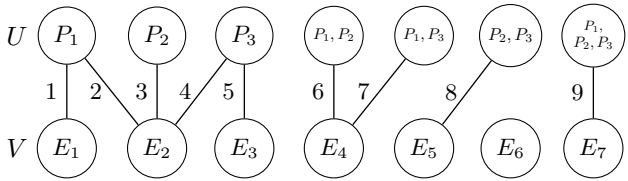
\begin{figure}[http]
	\centering
	\begin{tikzpicture}[scale=1, transform shape]
	\tikzstyle{edge_style} = [draw=black, line width=0.5]
	\node[text width=1cm] at (-0.5,0) {$U$};
	\node[text width=1cm] at (-0.5,-1.5) {{$V$}};
	\node[state,minimum size = 0.2cm] at (-0.2, 0) (p_nodeone) {\small{$P_1$}};
	\node[state,minimum size = 0.2cm] at (1.0, 0)     (p_nodetwo)     {\small$P_2$};
	\node[state,minimum size = 0.2cm] at (2.16, 0)     (p_nodethree)     {\small$P_3$};
	\node[state,scale = 0.65] at (3.35, 0)     (p_nodefour)     {\small$P_1,P_2$};
	\node[state,scale = 0.65] at (4.55, 0)     (p_nodefive)     {\small$P_1,P_3$};
	\node[state,scale = 0.65] at (5.7, 0) (p_nodesix) {$P_2,P_3$};
	\node[state,scale = 0.57, align=center] at (6.8, 0)     (p_nodeseven)     {\small$P_1,$\\$P_2,P_3$};
	\node[state,minimum size = 0.2cm] at (-0.2, -1.5)     (e_nodeone)     {\small$E_1$};
	\node[state,minimum size = 0.2cm] at (1.0, -1.5)     (e_nodetwo)     {\small$E_2$};
	\node[state,minimum size = 0.2cm] at (2.16, -1.5)     (e_nodethree)
	{\small$E_3$};
	\node[state,minimum size = 0.2cm] at (3.35, -1.5)     (e_nodefour)     {\small$E_4$};
	\node[state,minimum size = 0.2cm] at (4.55, -1.5)     (e_nodefive)     {\small$E_5$};
	\node[state,minimum size = 0.2cm] at (5.7, -1.5) (e_nodesix) {\small$E_6$};
	\node[state,minimum size = 0.2cm] at (6.8, -1.5) (e_nodeseven) {\small$E_7$};
	\draw[edge_style]  (p_nodeone) edge node[left] {\footnotesize{1}} (e_nodeone);
	\draw[edge_style]  (p_nodeone) edge
	node[left] {\footnotesize{2}}
	(e_nodetwo);
	\draw[edge_style]  (p_nodetwo) edge
	node[left] {\footnotesize{3}}
	(e_nodetwo);
	\draw[edge_style]  (p_nodethree) edge
	node[left] {\footnotesize{4}}
	(e_nodetwo);
	\draw[edge_style]  (p_nodethree) edge
	node[left] {\footnotesize{}5}
	(e_nodethree);
	\draw[edge_style]  (p_nodefour) edge
	node[left] {\footnotesize{6}}
	(e_nodefour);
	\draw[edge_style]  (p_nodefive) edge
	node[left] {\footnotesize{7}}
	(e_nodefour);
	\draw[edge_style]  (p_nodesix) edge
	node[left] {\footnotesize{8}}
	(e_nodefive);
	\draw[edge_style]  (p_nodeseven) edge
	node[left] {\footnotesize{9}}
	(e_nodeseven);
	\end{tikzpicture}
	\caption{The bipartite graph $G=(U\cup V,E)$ with 3 pursuers and 7 evaders, where the vertices in $U$ containing at least one common pursuer are conflicting in the maximum matching.}\label{MBMB_figure}
\end{figure}

\quad Let ${G}=({U}\cup {V},{E})$ be an undirected bipartite graph consisting of two independent vertex sets $U$ and $V$, where $E$ is the set of edges. We denote the edge connecting vertex $P_s\in U$ and vertex $E_j\in V$ by $e_{sj}$. In our problem, the vertex set $U$ consists of all pursuit coalitions of size no more than three, and $V$ represents the set of evaders. The bipartite graph $G$ is formally defined as follows:
\begin{equation}\begin{aligned}\label{origraphequ}
&U=[\mathcal{P}]^3,\quad V=\mathcal{E},\\
&E=\big\{e_{sj}\,|\,P_s\in U,E_j\in V, |\bar{\mathbb{E}}(s,j)\cap\Omega_{\rm goal}|\leq1,\\
&\qquad\qquad\qquad\quad\forall s_1\subsetneq s, |\bar{\mathbb{E}}(s_1,j)\cap\Omega_{\rm goal}|>1\big\}.
\end{aligned}\end{equation}
Notice that $|U|=N_o$, $|V|=N_e$. An edge $e_{sj}\in E$ if and only if pursuit coalition $P_s$ is able to capture $E_j$ in $\Omega_{\rm play}$ or at $\mathcal{T}$, while any subcoalition $s_1$ of $s$ cannot. An example of 3 pursuers and 7 evaders is depicted in Fig. \ref{MBMB_figure}.

\quad We aim to find a matching in the bipartite graph $G$ that contains a maximum number of evaders. However, since each pursuer can only appear in at most one pursuit coalition, the pursuit coalitions containing at least one common pursuer cannot coexist in the matching. As a result, our problem becomes a constrained maximum bipartite matching problem. We can also interpret the problem as an assignment problem with $N_p$ workers and $N_e$ jobs. In this assignment problem, some jobs are easy in the sense that they each can be finished by one individual worker, and some jobs are hard in the sense that they require cooperation among multiple workers\footnote{Even though we consider at most three workers (pursuers) for each job (evader), the derived results could be extended to the cases with any given number of workers for one job rather routinely.}. The goal is to find an assignment of workers to jobs such that as many jobs as possible are finished.

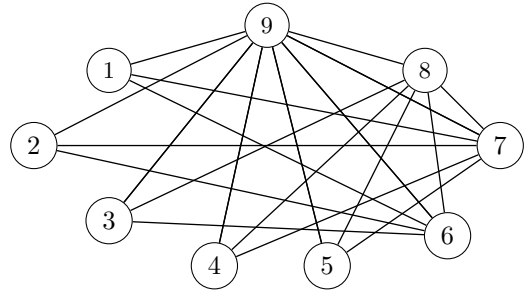
\begin{figure}[http]
	\centering
	\begin{tikzpicture}[scale=1, transform shape]
	\tikzstyle{edge_style} = [draw=black, line width=0.5]
	\node[state,minimum size = 0.2cm] at (3.5, 0) (node9)     {\small$9$};
	\node[state,minimum size = 0.2cm] at
	(1.4, -0.6) (node1) {\small$1$};
	\node[state,minimum size = 0.2cm] at
	(5.6, -0.6) (node8) {\small$8$};
	\node[state,minimum size = 0.2cm] at
	(0.4, -1.6) (node2) {$2$};
	\node[state,minimum size = 0.2cm] at
	(6.6, -1.6) (node7) {$7$};
	\node[state,minimum size = 0.2cm] at
	(1.4, -2.6) (node3) {$3$};
	\node[state,minimum size = 0.2cm] at
	(5.9, -2.8) (node6) {$6$};
	\node[state,minimum size = 0.2cm] at
	(2.8, -3.2) (node4) {$4$};
	\node[state,minimum size = 0.2cm] at
	(4.3, -3.2) (node5) {$5$};
	\draw[edge_style]  (node1) edge  (node6);
	\draw[edge_style]  (node1) edge  (node7);
	\draw[edge_style]  (node1) edge  (node9);
	\draw[edge_style]  (node2) edge  (node6);
	\draw[edge_style]  (node2) edge  (node7);
	\draw[edge_style]  (node2) edge  (node9);
	\draw[edge_style]  (node3) edge  (node6);
	\draw[edge_style]  (node3) edge  (node8);
	\draw[edge_style]  (node3) edge  (node9);
	\draw[edge_style]  (node3) edge  (node9);
	\draw[edge_style]  (node4) edge  (node7);
	\draw[edge_style]  (node4) edge  (node8);
	\draw[edge_style]  (node4) edge  (node9);
	\draw[edge_style]  (node4) edge  (node9);
	\draw[edge_style]  (node5) edge  (node7);
	\draw[edge_style]  (node5) edge  (node8);
	\draw[edge_style]  (node5) edge  (node9);
	\draw[edge_style]  (node5) edge  (node9);
	\draw[edge_style]  (node6) edge  (node8);
	\draw[edge_style]  (node6) edge  (node9);
	\draw[edge_style]  (node6) edge  (node9);
	\draw[edge_style]  (node7) edge  (node8);
	\draw[edge_style]  (node7) edge  (node9);
	\draw[edge_style]  (node7) edge  (node9);
	\draw[edge_style]  (node8) edge  (node9);
	\end{tikzpicture}
	\caption{The conflict graph $C$, where each vertex corresponds uniquely an edge of $G$. An edge in $C$ implies that its two related vertices (two edges in $G$) cannot coexist in the maximum matching of $G$.}
	\label{conflict_figure}
\end{figure}

\quad The conflicts among the pursuit coalitions can be represented by a conflict graph $C=(E,\bar{E})$ as in \cite{AI-MR-SLT:78,DJT:16,AD-UP-JS-GJW:11,TO-RZ-APP:13}. Each vertex in $C$ corresponds uniquely to an edge $e\in E$ of $G$. An edge $\bar{e}\in\bar{E}$ implies that the two vertices connected by $\bar{e}$ (two edges in $G$) cannot coexist in the maximum matching of $G$. The conflict graph $C$ may contain isolated vertices, which means that the corresponding edges in $G$ do not conflict with others. In our case, edges of $G$ incident to the vertices with at least one common pursuer are conflicting,
and thus the conflict graph $C$ is
\begin{equation}\label{conflictequ}
\bar{E}=\big\{(e_{sj},e_{pq})\,|\,e_{sj}\in E,e_{pq}\in E,s\neq p,s\cap p\neq\emptyset\big\}.
\end{equation}
For a better understanding, the conflict graph $C$ for our former example graph $G$ is depicted in Fig. \ref{conflict_figure}.

\quad Given the bipartite graph \eqref{origraphequ} and the conflict graph \eqref{conflictequ}, we define the binary integer programming (BIP) formulation for the maximum bipartite matching with conflict graph (MBMC) as follows:
\begin{equation}\label{BIP}
\begin{aligned}
& {\textup{maximize}}
& & \sum_{e_{sj}\in E}a_{sj} \\
& \textup{subject to}
& & \sum_{s\in U}a_{sj}\leq1,&&\forall E_j\in V,\\
&&& \sum_{j\in V}a_{sj}\leq1,&&\forall P_s\in U,\\
&&&a_{sj}+a_{pq}\leq1,&&\forall (e_{sj},e_{pq})\in\bar{E},\\
&&& a_{sj}\in\{0,1\}, &&\forall e_{sj}\in E,\\
&& &a_{sj}=0, &&\forall e_{sj}\notin E.
\end{aligned}
\end{equation}
where $a_{sj}=1$ indicates the assignment of pursuit coalition $P_s$ to capture $E_j$, and $a_{sj}=0$ means no assignment correspondingly. 

\quad Next, we prove the complexity of MBMC.

\begin{thom}\label{hardnessthom}{\rm (Hardness of the Matching).}
	The MBMC problem \eqref{BIP} is NP-hard.
\end{thom}
 \emph{Proof.} \rm
	We polynomially reduce the well-known NP-complete \emph{3-dimensional matching problem} \cite{RMK:72} to special instances of the MBMC problem. Let $\mathcal{I}=(X,Y,Z,T)$ be an arbitrary instance of 3-dimensional matching, where $X$, $Y$ and $Z$ are finite, disjoint sets with $|X|=|Y|=|Z|=m$, and $T$, a subset of $X\times Y\times Z$, consists of triples $(i,j,k)$ such that $i\in X$, $j\in Y$, and $k\in Z$. The problem is to determine whether there is a set $M\subseteq T$ such that $|M|=m$ and no two elements of $M$ agree in any coordinate. If so, the set $M$ is called a 3-dimensional matching of $\mathcal{I}$. We define the bipartite graph $G=(U\cup V, E)$ as follows:
	\[\begin{aligned}
	U=X\times Y,\  V=Z,\ E=\big\{((i,j),k)\,|\,(i,j,k)\in T\big\}.
	\end{aligned}\]
	
	\quad Let $M$ be a complete matching of $G$. If $M$ does not contain two edges $((i,j_1),k_1)$, $((i,j_2),k_2)$ or $((i_1,j),k_1)$, $((i_2,j),k_2)$, then $M$ corresponds to a 3-dimensional matching of $\mathcal{I}$. Therefore, the following restrictions are imposed:
	\[\begin{aligned}
	\bar{E}=\big\{\big(((i_1,j_1),k_1),((i_2,j_2),k_2)\big)\,|\,i_1=i_2\text{ or }j_1=j_2,\\
	((i_1,j_1),k_1)\in E,((i_2,j_2),k_2)\in E\big\}.
	\end{aligned}\]
	
	\quad Let the conflict graph $C=(E,\bar{E})$. Then, the matching problem in graph $G$ with conflict graph $C$ can be interpreted as a matching problem for the game with $2m$ (i.e., $|X|+|Y|$) pursuers and $m$ (i.e., $|Z|$) evaders. Each evader $k\in Z$ can be captured by two cooperative pursuers $i\in X$ and $j\in Y$ if $(i,j,k)\in T$. The pursuit coalitions with one or three pursuers do not exist in this case.
	
	\quad We now show that the MBMC problem $G$ with conflict graph $C$ has a complete matching if and only if $\mathcal{I}$ has a 3-dimensional matching.
	\begin{enumerate}
		\item Assume that $G$ with conflict graph $C$ has a complete matching $M$. Then, $|M|=m$ and $M$ is a subset of $E$. The conflict graph $C$ ensures that no two elements of $M$ agree in $X$ or $Y$ coordinate. Since $M$ is a matching, no two elements of $M$ agree in $Z$ coordinate. Therefore, $M$ is a 3-dimensional matching of $\mathcal{I}$ when we write $((i,j),k)$ as $(i,j,k)$.
		\item Let $M$ be a 3-dimensional matching of $\mathcal{I}$. For all $(i,j,k)\in M$, by definition, the edges $((i,j),k)$ constitute a complete matching for graph $G$ with conflict graph $C$.
	\end{enumerate}
	
	\quad Note that the decision problem of whether $G$ with conflict graph $C$ has a complete matching can be solved by computing the maximum matching of $G$ with conflict graph $C$. Thus, the 3-dimensional matching of $\mathcal{I}$ can be polynomially reduced to the MBMC problem, and the MBMC problem is NP-hard. \qed


\begin{algorithm}\label{SMA}
	\caption{Sequential Matching Algorithm}
	\begin{algorithmic}
		\State \textbf{Input:} A bipartite graph $G=(U\cup V, E)$ with $U=$ $[\mathcal{P}]^3$ and $V=\mathcal{E}$, where $e_{sj}\in E$ if the pursuit coalition $P_s$ in $U$ can defeat the evader $E_j$ in $V$
		\State \textbf{Output:}{ An approximation matching $M$ in $G$}
	\end{algorithmic}
	\begin{algorithmic}[1]
	\item $U_1\leftarrow\mathcal{P}$, $V_1\leftarrow\mathcal{E}$, $E_1\leftarrow\big\{e_{sj}\in E\,|\,P_s\in U_1,E_j\in V_1\big\}$
	\item Compute the maximum matching $M_1$ in the subgraph $G_1=(U_1\cup V_1,E_1)$ by maximum network flow \cite{LRF-DRF:62}; 
	\item $A_1\leftarrow\{P_i\,|\,i\in s,e_{sj}\in M_1\},B_1\leftarrow\{E_j\,|\,e_{sj}\in M_1\}$ 
	\item $U_2\leftarrow[\mathcal{P}\setminus A_1]^2$,
	$V_2\leftarrow\mathcal{E}\setminus B_1$ 
	\item $E_2\leftarrow\{e_{sj}\in E\,|\,P_s\in U_2,E_j\in V_2\}$
	\item Compute the maximum matching $M_2$ in the subgraph $G_2=(U_2\cup V_2,E_2)$ by maximum network flow;
	\item $A_2\leftarrow\{P_i\,|\,i\in s,e_{sj}\in M_2\}, B_2\leftarrow\{E_j\,|\,e_{sj}\in M_2\}$
	\item $U_3\leftarrow\big[\mathcal{P}\setminus \big(A_1\cup A_2\big)\big]^3,V_3\leftarrow\mathcal{E}\setminus \big(B_1\cup B_2\big)$
	\item $E_3\leftarrow\{e_{sj}\in E\,|\,P_s\in U_3, E_j\in V_3\}$
	\item Compute the maximum matching $M_3$ in the subgraph $G_3=(U_3\cup V_3,E_3)$ by maximum network flow;
	\item \textbf{Return} $M=M_1\cup M_2\cup M_3$.
	\end{algorithmic}
\end{algorithm}

\quad Next, we give an approximation algorithm called Sequential Matching Algorithm stated in Algorithm \ref{SMA} for MBMC. We sketch out the main idea of the Sequential Matching Algorithm as follows:\\
Step 1 (from line 1 to 2): Use maximum network flow to compute the maximum matching $M_1$ of the subgraph $G_1$ which only considers the pursuit vertices containing one pursuer.\\
Step 2 (from line 3 to 6): Let $A_1$ and $B_1$ be the sets of pursuers and evaders in $M_1$ respectively. Remove the vertices of $G$ containing at least one player occurring in the set $A_1\cup B_1$, and for the remaining part, construct the subgraph $G_2$ which only considers the pursuit vertices containing at most two pursuers. Then, use maximum network flow to compute the maximum matching $M_2$ of $G_2$. \\
Step 3 (from line 7 to 10): Let $A_2$ and $B_2$ be the sets of pursuers and evaders in $M_2$ respectively. Remove the vertices of $G$ containing at least one player occurring in the set $A_1\cup A_2\cup B_1\cup B_2$, and for the remaining part, obtain the subgraph $G_3$. Then, use maximum network flow to compute the maximum matching $M_3$ of $G_3$.\\
Step 4: The output is the union of these three matchings.

\quad It turns out that this algorithm has  great features.

\begin{thom}\label{confactorthom}{\rm (Constant-Factor Approximation Algorithm).}
	The Sequential Matching Algorithm is of polynomial time and
	\begin{enumerate}[(i)]
		\item a 1/3-approximation algorithm for MBMC;
		\item a 1/2-approximation algorithm if the solution of MBMC does not contain pursuit coalitions with three pursuers;
		\item an exact algorithm if the solution of MBMC does not contain pursuit coalitions with two or three pursuers.
	\end{enumerate}
\end{thom}

 \emph{Proof.} \rm
	We postpone the proof to Appendix~\ref{appdix:approximationalgo}. \qed

\begin{cor}{\rm (Class of Complexity on Approximation Algorithm).}
	The MBMC is APX-complete.
\end{cor}
 \emph{Proof.}  \rm
	Note that 3-dimensional matching can be polynomially reduced to the MBMC. Since 3-dimensional matching is APX-complete \cite{VK:91} and \thomref{confactorthom} implies that the MBMC is in APX, the statement follows. \qed

\begin{algorithm}
	\caption{Receding Horizon Strategy\label{al2}}
	\ \ \  \textbf{Input:} $\mathcal{P},\mathcal{E},\{\mathbf{x}_{P_i}\}_{P_i\in\mathcal{P}},\{\mathbf{x}_{E_j}\}_{E_j\in\mathcal{E}},M_a\leftarrow\emptyset,L_c\leftarrow0$
	\begin{algorithmic}[1]
		\item \textbf{Repeat}
		\item \quad $M\leftarrow\text{SeqMax}\big(\mathcal{P},\mathcal{E},\{\mathbf{x}_{P_i}\}_{P_i\in\mathcal{P}},\{\mathbf{x}_{E_j}\}_{E_j\in\mathcal{E}}\big)$ 
		\item \quad\textbf{if} $|M|>|M_a|$\textup{ or }$L_c==1$ \textbf{then}
		\item \qquad $M_a\leftarrow M,L_c\leftarrow0$
		\item \quad \textbf{end if}
		\item \quad 
			Assign a pursuit coalition to each evader that is \text{ \ \ } part of the matching $M_a$;
	
		\item \quad For a short duration $\Delta$, apply the ES-based strat- \text{\ \ \, } egy for each pursuer that is part of the matching \text{ \ \ } $M_a$. For the rest of the pursuers and for all evaders \text{ \ \ } in $\mathcal{E}$, apply some (any) strategy;
		\item \quad Update the player positions after the duration $\Delta$;
		\item \quad\textbf{for} every evader $E_j $  in $\mathcal{E}$ \textbf{do}
		\item \qquad \textbf{if} $E_j$ is captured or enters $\Omega_{\rm goal}$ \textbf{then}
		\item \quad \qquad \textbf{if} $E_j$ is captured \textbf{then}
		\item \qquad\qquad $L_c\leftarrow1$
		\item \quad \qquad \textbf{end if}
		\item \quad \qquad $\mathcal{E}\leftarrow\mathcal{E}\setminus E_j$
		\item \qquad \textbf{end if}
		\item \quad \textbf{end for}
		\item \textbf{until} $\mathcal{E}=\emptyset$;
	\end{algorithmic}
\end{algorithm}

%

\subsection{Receding Horizon Strategy}
\quad In this subsection, we design a receding horizon strategy for the pursuit team. This strategy is useful because a better matching  may occur as the game runs, and a rematching should be performed when an evader is captured. With Algorithm \ref{al2}, the bipartite graph and the corresponding approximate maximum matching can be updated, potentially in real time, as players change their positions during the game. $M_a$ and $L_c$ denote the adopted matching and label for the capture of evaders, respectively, and $\text{SeqMax}\big(\mathcal{P},\mathcal{E},\{\mathbf{x}_{P_i}\}_{P_i\in\mathcal{P}},\{\mathbf{x}_{E_j}\}_{E_j\in\mathcal{E}}\big)$ computes the maximum matching by Sequential Matching Algorithm.

\quad This paper considers a constant time step $\Delta>0$. In general, the solution to the maximum matching is not unique. Besides,
as the players play out the game in real time, more alternative solutions may generate. Thus, to avoid unnecessary change of the interception scheme, assume that the pursuit team would never change the former matching unless more evaders can be captured in the new
matching as lines 3-5 show. As long as each pursuer uses the ES-based strategy from the related pursuit coalition against the matched evader, the size of the matching never decreases until an evader is captured.

\section{Bounded Convex Play Region with An Exit}\label{BoundedSec}
\begin{figure}
	\centering
	\graphicspath{{figure_original/}}
	\subfigure{
		\includegraphics[width=82mm,height=32mm]{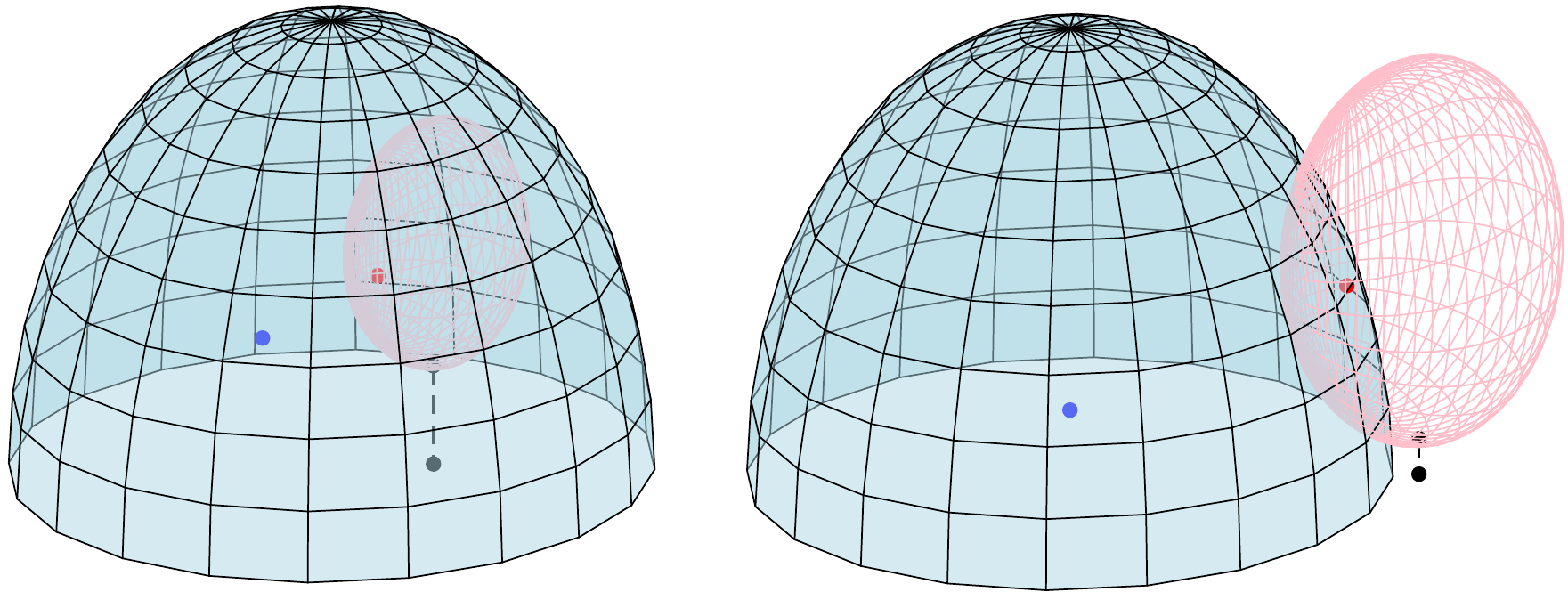}
		\put(-199,6){\scriptsize{$\mathcal{T}^b$}}
		\put(-90,5){\scriptsize{$\mathcal{T}^b$}}
		\put(-233,79){\scriptsize{$\Omega_{\rm play}^b$}}
		\put(-129,72){\scriptsize{$\Omega_{\rm play}^b$}}
		\put(-176,58){\scriptsize{$\bar{\mathbb{E}}(i,j)$}}
		\put(-75,-6){\scriptsize{$(b)$}}
		\put(-192,-5){\scriptsize{$(a)$}}
		\put(-32,85){\scriptsize{$\bar{\mathbb{E}}(i,j)$}}
	}
	\caption{Bounded convex play region with a planar exit, where the interception point $I^b(i,j)$ lies $(a)$ in the interior of the play region; $(b)$ at the boundary of the play region.}
	\label{bounded_figure}
\end{figure}
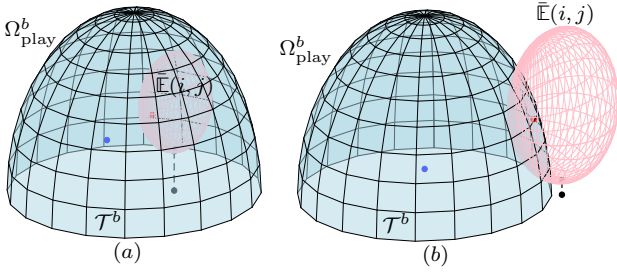

\quad In this section, we extend the previous analysis to the case when the game is played in a 3D bounded convex region. We consider a bounded convex play region with an exit through which the evaders escape from the play region. The exit is assumed to be a part of a plane. The goal of the pursuit team is to capture as many evaders as possible before the evaders leave the play region. The play region $\Omega_{\rm play}^b$ is a closed convex region in $\mathbb{R}^3$ and the exit $\mathcal{T}^b$ is a part of its boundary. Formally,
\begin{equation}\begin{aligned}
\Omega_{\rm play}^b&=\big\{\mathbf{x}\in\mathbb{R}^3\,|\,z\ge0,g(\mathbf{x})\ge0\big\},\\
\mathcal{T}^b&=\big\{\mathbf{x}\in\Omega_{\rm play}^b\,|\,z=0\big\},
\end{aligned}\end{equation}
where $g:\mathbb{R}^3\rightarrow\mathbb{R}$ is a differentiable function such that the set $\Omega_1^b=\{\mathbf{x}\in\mathbb{R}^3\,|\,g(\mathbf{x})\ge0\}$ is convex. Additionally, Assume that $\Omega_{\rm play}^b$ is non-empty and contains more than one point. Condition 4) in \aspref{IsolationDeploy} becomes that $\mathbf{x}^0_{P_i}\in\Omega_{\rm play}^b$ for all $P_i\in\mathcal{P}$ and $\mathbf{x}^0_{E_j}\in\Omega_{\rm play}^b$ for all $E_j\in\mathcal{E}$. An example of the play region is given in Fig. \ref{bounded_figure}.

\quad For $P_s\in[\mathcal{P}]^{+}$ and $E_j\in\mathcal{E}$, if $\bar{\mathbb{E}}(s,j)\cap\mathcal{T}^b=\emptyset$, then we define the interception point $I^b(s,j)\in\bar{\mathbb{E}}(s,j)\cap\Omega_{\rm play}^b$ as the closest point to the plane containing $\mathcal{T}^b$, and $I^b(s,j)$ is the solution to the following convex problem
\begin{equation}\label{convexpbm2}
\begin{aligned}
& \underset{\mathbf{x}\in\mathbb{R}^3}{\textup{minimize}}
& & z \\
& \textup{subject to}
& & f_{ij}(\mathbf{x})\ge0,\quad \forall i\in s,\\
& & & g(\mathbf{x})\ge0,
\end{aligned}
\end{equation}
where the constraint $z\ge0$ in $\Omega^b_{\rm play}$ is not involved because it holds naturally when $\bar{\mathbb{E}}(s,j)\cap\mathcal{T}^b=\emptyset$.

\begin{lema}\label{uniquegoalpoint2}{\rm (Uniqueness of the Interception Point).}
	For any $P_s\in[\mathcal{P}]^{+}$ and $E_j\in\mathcal{E}$, if $\bar{\mathbb{E}}(s,j)\cap\mathcal{T}^b=\emptyset$, the interception point $I^b(s,j)$ is unique.
\end{lema}
 \emph{Proof.} \rm
	Since $\bar{\mathbb{E}}(s,j)\cap\mathcal{T}^b=\emptyset$, similar to \defiref{goalpointlema}, we denote $I(s,j)$ as the unique point in $\bar{\mathbb{E}}(s,j)$ that is closest to the plane containing $\mathcal{T}^b$.
	If $I(s,j)$ lies in $\Omega_{\rm play}^b$ as Fig. \ref{bounded_figure}(a) shows, then $I^b(s,j)=I(s,j)$ and thus $I^b(s,j)$ is unique. If $I(s,j)$ lies out of $\Omega_{\rm play}^b$ as Fig. \ref{bounded_figure}(b) indicates, we consider a plane $\mathcal{T}_1$ parallel to $\mathcal{T}^b$ and move it from $\mathcal{T}^b$ towards $\Omega_{\rm play}^b$. At the beginning, $\mathcal{T}_1$'s intersection sets with $\bar{\mathbb{E}}(s,j)$ and $\Omega_{\rm play}^b$ are two disjoint sets: a strictly convex set and a convex set respectively. As $\mathcal{T}_1$ moves, these two intersection sets are tangent when they intersect at the first time. The tangent point is $I^b(s,j)$ and the statement follows from the uniqueness of $I^b(s,j)$. \qed

\quad In the next, we can still construct a similar ES-based strategy to guarantee the winning of the pursuer coalitions.

\begin{lema}\label{ES-based2}{\rm (ES-Based Strategy).}
	For any $P_s\in[\mathcal{P}]^{+}$ and $E_j\in\mathcal{E}$, suppose that $\bar{\mathbb{E}}(s,j)\cap\mathcal{T}^b=\emptyset$. If every pursuer $P_i$ in $P_s$ adopts the feedback strategy $\mathbf{u}_{P_i}=\frac{I^b(s,j)-\mathbf{x}_{P_i}}{\|I^b(s,j)-\mathbf{x}_{P_i}\|_2}$, then it is guaranteed that $\bar{\mathbb{E}}(s,j)$ does not approach $\mathcal{T}^b$, i.e., $\dot{z}_{I^b(s,j)}\ge0$ for any $\mathbf{u}_{E_j}\in\mathbb{S}^2$. Moreover, $\dot{z}_{I^b(s,j)}=0$ if and only if $E_j$ adopts the feedback strategy $\mathbf{u}_{E_j}=\frac{I^b(s,j)-\mathbf{x}_{E_j}}{\|I^b(s,j)-\mathbf{x}_{E_j}\|_2}$.
\end{lema}
 \emph{Proof.} \rm
	Since $\bar{\mathbb{E}}(s,j)\cap\mathcal{T}^b=\emptyset$, by \lemaref{uniquegoalpoint2}, the interception point $I^b(s,j)$ is unique. If $I(s,j)$ lies inside of $\Omega_{\rm play}^b$ as Fig. \ref{bounded_figure}(a) shows, then $I^b(s,j)=I(s,j)$ and the statement follows from \thomref{ESstramultilema}. In the following, we focus on the case when $I(s,j)$ lies outside of $\Omega_{\rm play}^b$ as in Fig. \ref{bounded_figure}(b). For simplicity, we denote $I(s,j)$ and $I^b(s,j)$ by $\mathbf{x}_I$ and $\mathbf{x}_I^b$, respectively.
	
	\quad Since $\mathbf{x}_I^b$ moves on $\partial\Omega_1^b$, i.e., $g(\mathbf{x}_I^b)\equiv0$, thus we have
	\begin{equation}\label{boundequ1}
	\frac{dg(\mathbf{x}_I^b)}{dt}    =0\Rightarrow\nabla g(\mathbf{x}_I^b)^\top\dot{\mathbf{x}}_I^b=0.
	\end{equation}
	The KKT conditions for (\ref{convexpbm2})  are as follows:
	\begin{equation}\begin{aligned}\label{KKT2}
	&[
	0\ \ 0\ \ -1
	]^\top
	=\sum_{i\in s}\lambda_i\nabla f_{ij}(\mathbf{x}^b_I)+\lambda_g\nabla g(\mathbf{x}_I^b),\\
	&f_{ij}(\mathbf{x}_I^b)\ge0,\,\lambda_i\leq0,\,\lambda_if_{ij}(\mathbf{x}_I^b)=0,\quad \forall i\in s,\\
	& \lambda_g\leq0,\, g(\mathbf{x}_I^b)\ge0,\,\lambda_gg(\mathbf{x}_I^b)=0,
	\end{aligned}\end{equation}
	where $\lambda_g\in\mathbb{R}$ is the Lagrange multiplier related to $g(\mathbf{x})\ge0$. Additionally, similar to (\ref{oneESbasedinequ}), we obtain
	\begin{equation}\begin{aligned}\label{timederivativeequ22}
	&\nabla f_{ij}(\mathbf{x}^b_I)^\top\dot{\mathbf{x}}_I^b\ge0\text{ for all }i\in s \textup{ with }f_{ij}(\mathbf{x}_I^b)=0,\end{aligned}\end{equation}when $P_i$ adopts the feedback strategy $\mathbf{u}_{P_i}=\frac{\mathbf{x}_I^b-\mathbf{x}_{P_i}}{\|\mathbf{x}^b_I-\mathbf{x}_{P_i}\|_2}$.
	Thus, it follows from \eqref{slackcondition}, \eqref{boundequ1}, \eqref{KKT2} and \eqref{timederivativeequ22} that $z_I^b$ satisfies \[\begin{aligned}
	-\dot{z}^b_I&=
	[
	0\ \ 0\ \ -1
	]
	\dot{\mathbf{x}}_I^b\\
	&=\sum_{i\in s}\lambda_i\nabla f_{ij}(\mathbf{x}^b_I)^\top\dot{\mathbf{x}}_I^b+\lambda_g\nabla g(\mathbf{x}_I^b)^\top\dot{\mathbf{x}}_I^b\leq0,
	\end{aligned}\]
	which leads to the similar conclusion as \thomref{ESstramultilema}.\qed

\quad Then, the results about the game of kind are straightforward, formally stated below.

\begin{rek}
	Similarly, the game winner between a pursuit team $P_s$ and an evader $E_j$ for the bounded convex play region is determined as follows: If $\bar{\mathbb{E}}(s,j)\cap\mathcal{T}^b$ is empty, then the pursuit team $P_s$ wins; if $\bar{\mathbb{E}}(s,j)\cap\mathcal{T}^b$ has more than one element, then $E_j$ wins; if $\bar{\mathbb{E}}(s,j)\cap\mathcal{T}^b$ has a unique element, then two teams are tied.
\end{rek}


\quad Since $\bar{\mathbb{E}}(s,j)\cap\Omega_{\rm play}^b\subseteq\bar{\mathbb{E}}(s,j)$, similar to the case of unbounded play region, we need at most three pursuers to capture one evader before the latter escapes. Thus, the results of maximum matching can also be applied.

\begin{figure}
	\centering
	\graphicspath{{figure_original/}}
	\subfigure{
		\includegraphics[width=78mm,height=40mm]{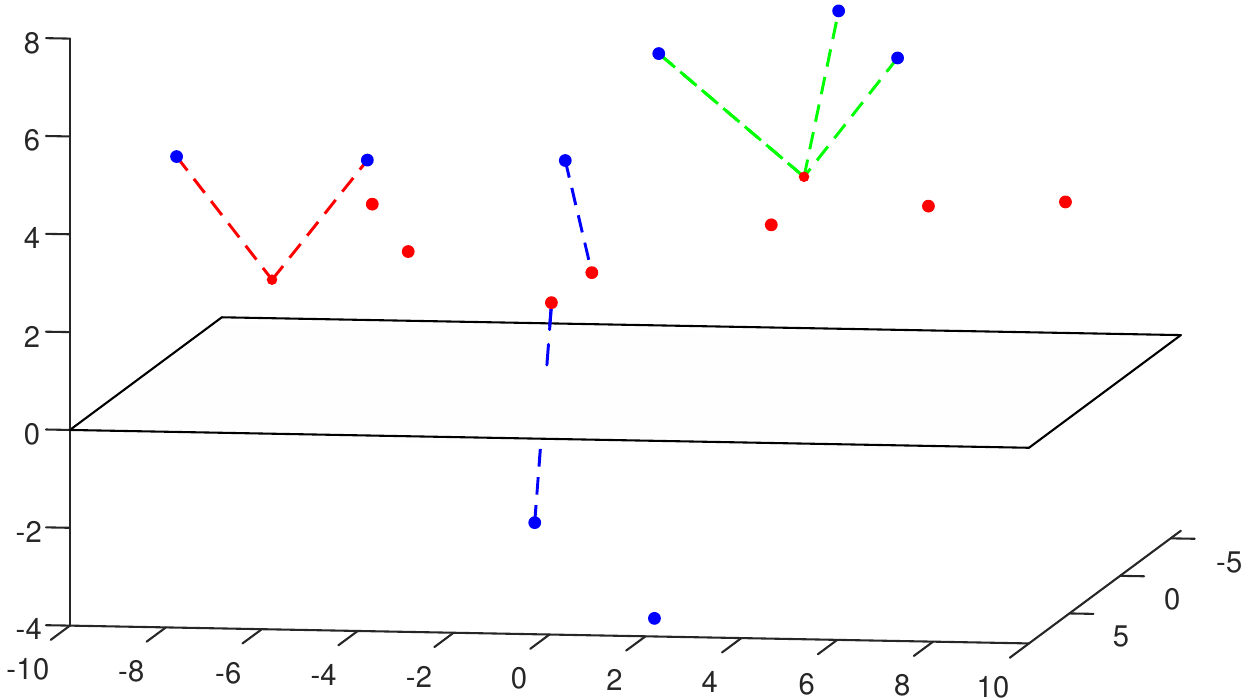}
		\put(-127,-5){\scriptsize{$(a)$}}
		\put(-206,105){\scriptsize{$z$}}
	}
	\subfigure{
		\includegraphics[width=78mm,height=40mm]{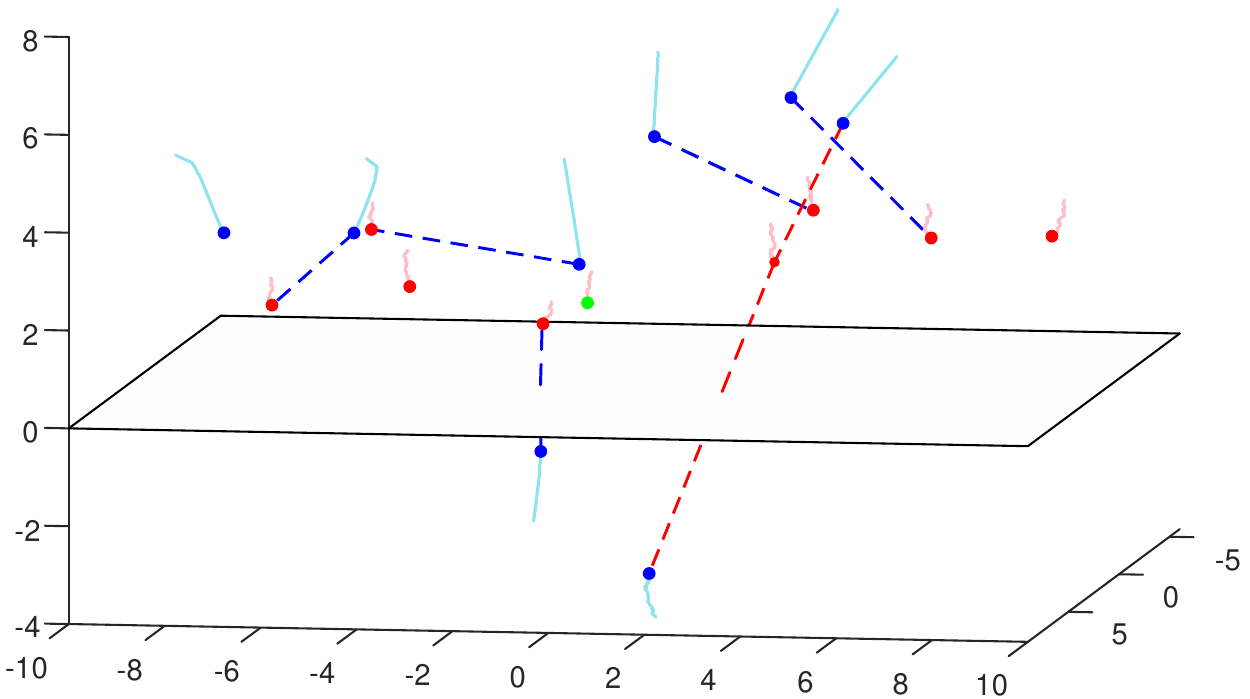}
		\put(-127,-5){\scriptsize{$(b)$}}
		\put(-206,105){\scriptsize{$z$}}
	}
	\subfigure{
		\includegraphics[width=78mm,height=40mm]{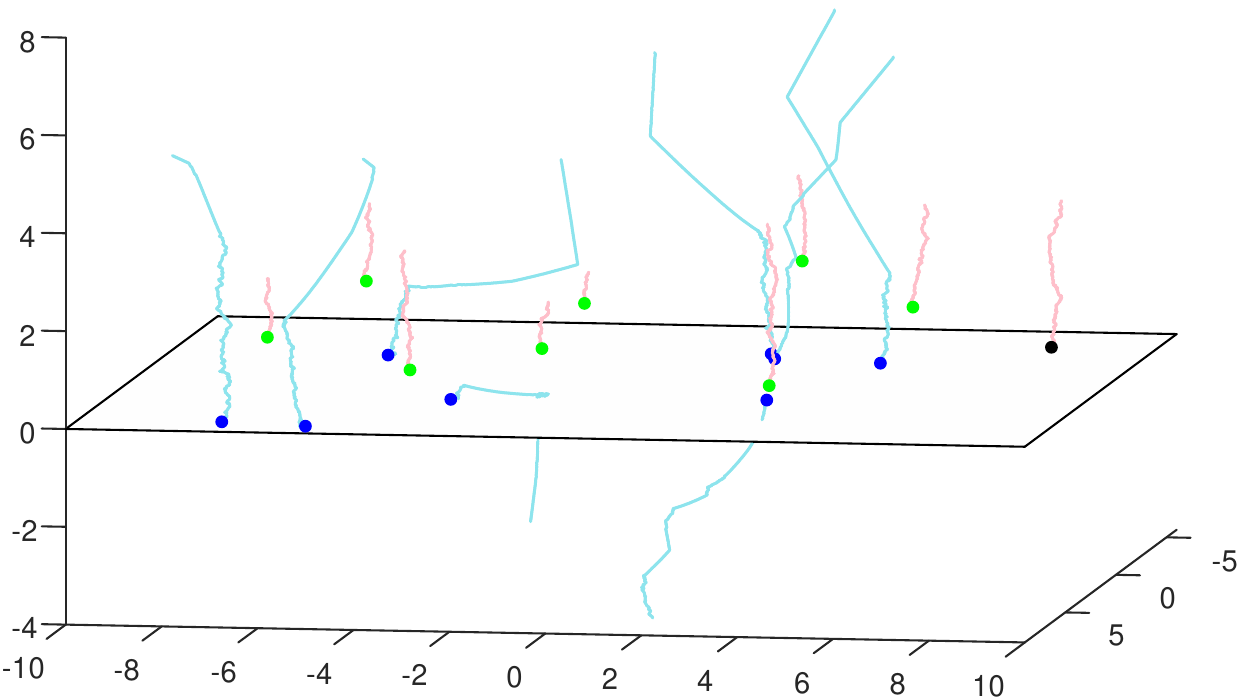}
		\put(-127,-5){\scriptsize{$(c)$}}
		\put(-206,105){\scriptsize{$z$}}
	}
	\subfigure{
		\includegraphics[width=82mm,height=14mm]{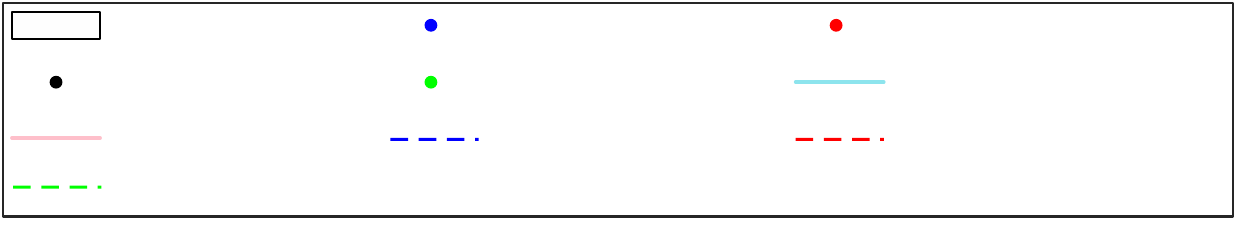}
		\put(-212,33){\tiny{Separating Plane}}
		\put(-212,23.5){\tiny{Arrived Evader}}
		\put(-212,14){\tiny{Evader Traj.}}
		\put(-212,4.5){\tiny{3-to-1 Matching}}
		\put(-141, 14){\tiny{1-to-1 Matching}}
		\put(-141,23.5){\tiny{Captured Evader}}
		\put(-141,33){\tiny{Pursuer}}
		\put(-63,14){\tiny{2-to-1 Matching}}
		\put(-64.5,33){\tiny{Uncaptured Evader}}
		\put(-63,23.5){\tiny{Pursuer Traj.}}
	}
	\caption{Simulation of a game with eight pursuers and nine evaders in the unbounded play region, showing that one evader succeeds to reach the goal region ($z\leq0$) and eight evaders are captured in the play region ($z>0$).}
	\label{bounded_example_fig}
\end{figure}

\section{Numerical Results}\label{simulationsec}
\quad In this section, numerical results are presented to illustrate the previous theoretical developments for the cases of unbounded and bounded convex play regions. Numerical studies are performed in Matlab R2017b on a laptop with a Core i7-8550U processor with 16 GB of memory.
\subsection{Unbounded Play Region}
\quad We first consider the unbounded play region with $N_p=8$ and $N_e=9$. Initially, as Fig. \ref{bounded_example_fig}(a) shows, the maximum-matching strategy indicates that four evaders are matched by seven pursuers, including two 1-to-1, one 2-to-1 and one 3-to-1 matchings. The matched evaders will be captured, unless the pursuit team changes its matching when a matching of greater size occurs as the game runs. A snapshot of the game is presented in Fig. \ref{bounded_example_fig}(b) where the pursuit team changes its matching because a better matching with six matched evaders occurs. In the end, as shown in Fig. \ref{bounded_example_fig}(c), one evader reaches the goal region successfully and eight evaders are captured in the play region.
\subsection{Bounded Convex Play Region}
\begin{figure}
	\centering
	\graphicspath{{figure_original/}}
	\subfigure{
		\includegraphics[width=80mm,height=47mm]{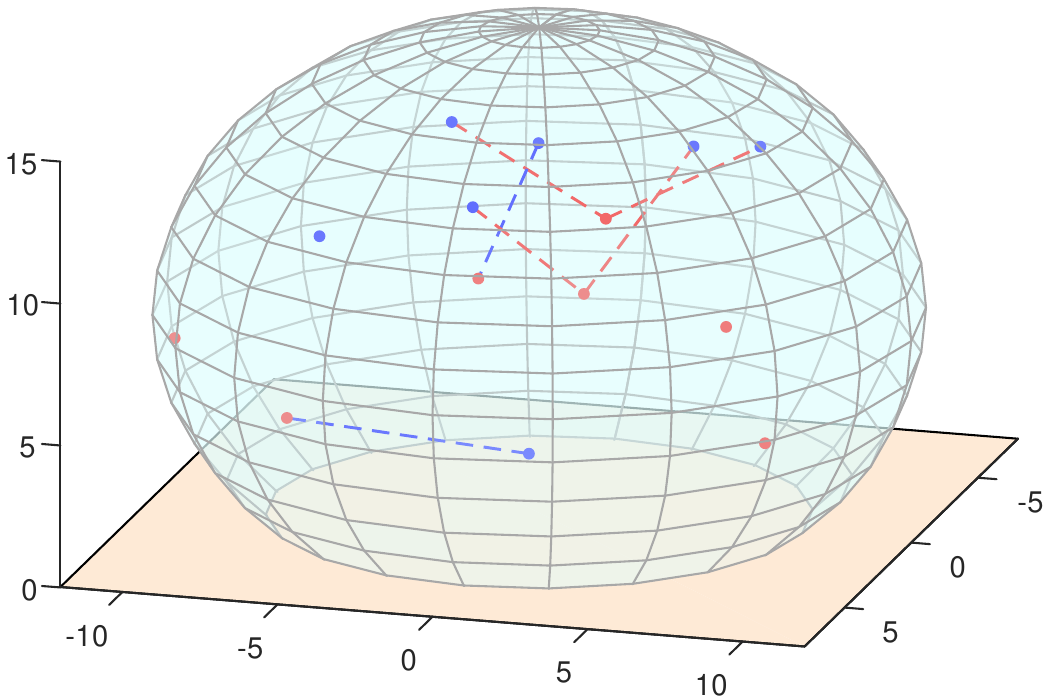}
		\put(-127,-2){\scriptsize{$(a)$}}
		\put(-211,99){\scriptsize{$z$}}
	}
	\subfigure{
		\includegraphics[width=80mm,height=47mm]{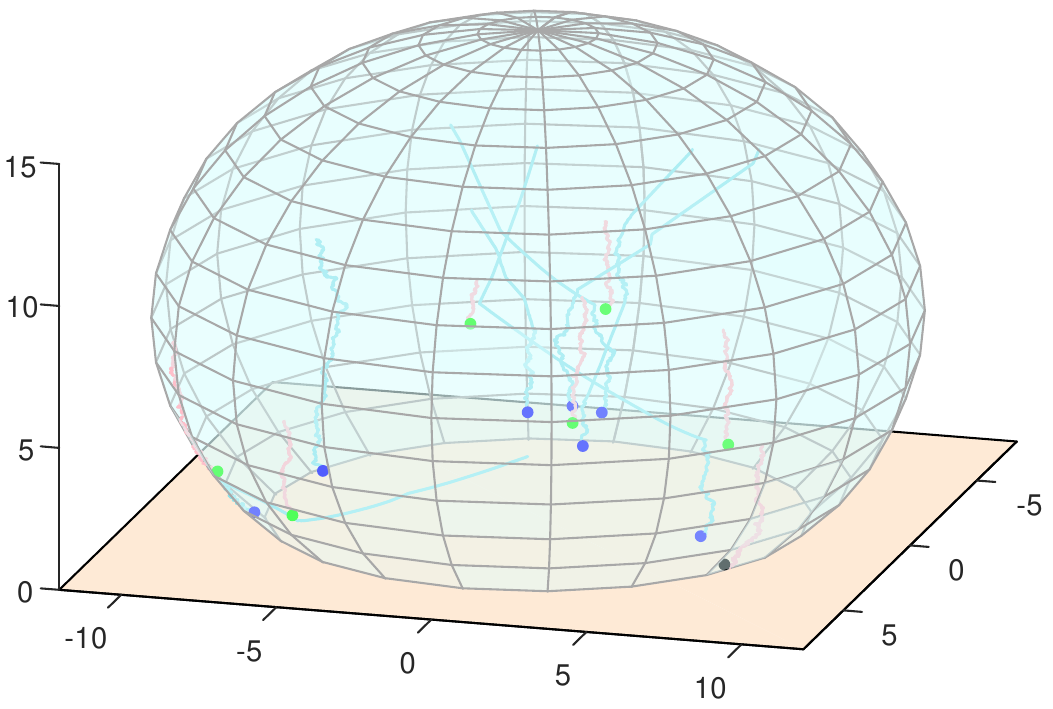}
		\put(-127,-2){\scriptsize{$(b)$}}
		\put(-211,99){\scriptsize{$z$}}
	}
	\subfigure{
		\includegraphics[width=83mm,height=5mm]{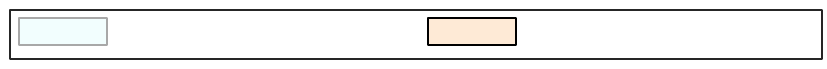}
		\put(-203,5){\tiny{Boundary of Play Region}}
		\put(-88.2,5){\tiny{Exit Plane of Play Region}}
	}
	\caption{Simulation of a game with seven pursuers and seven evaders in the bounded convex play region, showing that one evader succeeds to escape from the play region and six evaders are captured.}
	\label{bounded_example_fig2}
\end{figure}

\quad In this section, we consider a bounded convex play region with a planar exit, where $N_p=7$ and $N_e=7$. As Fig. \ref{bounded_example_fig2}(a) shows, four evaders are matched by six pursuers at the beginning, including two 1-to-1 and two 2-to-1 matchings. Finally, the pursuit team captures six evaders successfully and one evader escapes.

\section{Conclusion}\label{conclusionsec}
\quad In this paper, we studied a 3D multiplayer reach-avoid game where multiple pursuers defend a goal region against multiple evaders. For multiple pursuers and one evader, we showed that the evasion space corresponding to a pursuit coalition and an evader is strictly convex and the associated interception point is unique and can be computed via a convex program. We further revealed that the pursuit coalition can always defend the goal region by moving towards the interception point if the initial condition allows. We also found that in 3D if a pursuit coalition can defend the goal region against an evader,  then at most three pursuers in the coalition are necessarily needed. We solved the HJI equation associated with a special subgame of degree by a convex program. For multiple pursuers and multiple evaders, the matching is considered. We have shown that our matching between pursuit coalitions and evaders is an instance of a class of constrained matching problems, i.e., MBMC. We analyzed the complexity of MBMC and designed a constant-factor approximation algorithm with polynomial computation time to solve it. We also demonstrated that our results can be applied to the case of a bounded convex play region by slightly modifying the interception point. Future work will focus on distributed multiplayer reach-avoid games.

\begin{appendices}
	\section{Proof of \lemaref{polarcoor}}\label{appx:convexcurve}
	\quad The curvature $\kappa$ of the curve $\rho=\rho(\psi)$ in polar coordinates is given by \cite[Lemma 3.7]{EA-SS-AG:17}
	\begin{equation*}
	\kappa=\frac{\rho^2+2(\frac{d\rho}{d\psi})^2-\rho\frac{d^2\rho}{d\psi^2}}{\big(\rho^2+(\frac{d\rho}{d\psi})^2\big)^{3/2}}.
	\end{equation*}
	Therefore, if \eqref{eq:curvature} holds, then we have $\kappa>0$ for all $\psi\in[0,2\pi]$, and by \cite[Problem 1.7.6]{VAT:06}
	, the curve is convex. Moreover, since $\kappa$ is strictly positive, the curve does not contain any line segments and thus is strictly convex.
	
	\quad Since the curve is convex, by definition, the set consisting of the curve and its interior has a supporting hyperplane at every point on the boundary. Along with the fact that the set is closed and has nonempty interior, we have that the set is convex. The strict convexity of the set follows from that of the curve.
	
	\section{Proof of \thomref{confactorthom}}\label{appdix:approximationalgo}
	\quad Let $G=(U\cup V, E)$ with conflict graph $C$ be an instance of MBMC, where $G$ and $C$ are given by (\ref{origraphequ}) and (\ref{conflictequ}), respectively.
	
	\quad Assume that the Sequential Matching Algorithm is applied on this instance, and returns a matching $M=M_1\cup M_2\cup M_3$, where $M_1$, $M_2$ and $M_3$ may be empty. As lines 4 and 8 in Algorithm \ref{SMA} show, we remove the matched pursuers and evaders when computing the next matching. Thus, $M$ satisfies the conflict graph $C$ naturally.
	
	\quad Since this algorithm involves solving three maximum network flow problems which can be solvable in polynomial time \cite{LRF-DRF:62}, it is also of polynomial time.
	
	\quad Let $M^*=M_1^*\cup M_2^*\cup M_3^*$ be an optimal solution of MBMC on the given instance, where $M_i^*$ ($i=1,2,3$) is the set of edges incident to vertices in $U$ with $i$ pursuers. Next, we give an upper bound of $|M_i^*|$.
	
	\quad $(i)$ Since $M_1$ is the maximum matching of the subgraph $G_1$, we have
	\begin{equation}\label{approxiequ1}
	|M_1^*|\leq |M_1|=|A_1|=|B_1|\leq|M|.
	\end{equation}
	
	\quad $(ii)$ For the subgraph $G_2$, we add a pursuer $P_i\in A_1$ into $U_2$, and obtain a new subgraph $G_2'=(U_2'\cup V_2', E_2')$ of $G$ with
	\[\begin{aligned}
	&U_2'=[\mathcal{P}/A_1\cup P_i]^2,\quad V_2'=V_2,\\ &E_2'=\big\{e_{sj}\in E\,|\,P_s\in U_2',E_j\in V_2'\big\}.
	\end{aligned}
	\]
	From now on, we omit the formulation of the conflict graph for any graph we construct, because it can be obtained routinely by at most one appearance of each pursuer in the matching.
	
	\quad Let $M_2'$ be the maximum matching of $G_2'$, and it is easy to see that $|M_2'|\ge|M_2|$. Note that pursuer $P_i$ can capture at most one evader by itself or cooperation with another pursuer. If we remove pursuer $P_i$ from $U_2'$, $|M_2'|$ decreases by at most 1 and $G'_2$ is reduced to $G_2$. Thus, we have $|M_2'|\leq|M_2|+1$.
	
	\quad The similar results can be obtained when we add an evader $E_j\in B_1$ into $V_2$. Thus, by adding all pursuers in $A_1$ and all evaders in $B_1$ into the graph $G_2$, we can obtain the subgraph $G''_2=(U''_2\cup V''_2, E''_2)$ of $G$ with
	\[
	U''_2=[\mathcal{P}]^2,\ V''_2=\mathcal{E},\ E_2''=\big\{e_{sj}\in E\,|\,P_s\in U_2'',E_j\in V_2''\big\},
	\]
	and the maximum matching $M''_2$ of $G_2''$ satisfies
	\[
	|M_2''|\leq|M_2|+|A_1|+|B_1|.
	\]
	The graph $G_2''$ consists of all edges of $G$ incident to vertices in $U$ containing one or two pursuers. Thus, $|M_1^*|+| M_2^*|$ is bounded by
	\begin{equation}\label{approxiequ2}
	|M_1^*|+| M_2^*|\leq|M_2''|\leq|M_2|+2|M_1|\leq2|M|.
	\end{equation}
	
	\quad $(iii)$ For the subgraph $G_3$, we add a pursuer $P_i\in A_1\cup A_2$ into $U_3$ and obtain a new subgraph
	$G_3'=(U_3'\cup V_3', E_3')$ of $G$ with
	\[\begin{aligned}
	&U_3'=\big[\mathcal{P}\setminus \big(A_1\cup A_2\big)\cup P_i\big]^3,\quad V_3'=V_3,\\ &E_3'=\big\{e_{sj}\in E\,|\,P_s\in U_3',E_j\in V_3'\big\}.
	\end{aligned}
	\]
	Analogously, the maximum matching $M_3'$ of $G_3'$ satisfies $|M_3'|$ $\leq|M_3|+1$. By putting all pursuers in $A_1\cup A_2$ and all evaders in $B_1\cup B_2$ into $G_3$, then $G_3$ becomes $G$ and its maximum matching $M^*$ satisfies
	\begin{equation}\begin{aligned}\label{approxiequ3}
	|M^*|&\leq|M_3|+|A_1|+|A_2|+|B_1|+|B_2|\\
	&=|M_3|+2|M_1|+3|M_2|\leq3|M|,
	\end{aligned}
	\end{equation}where we have used $|A_1|=|B_1|=|M_1|$ and  $|A_2|=2|B_2|=2|M_2|$.
	
	Thus, the statement follows from \eqref{approxiequ1}, \eqref{approxiequ2} and \eqref{approxiequ3}.

\end{appendices}


 \small
 \bibliographystyle{plain}
\bibliography{bib/alias,bib/Main,bib/FB,bib/New}          


\end{document}